\documentclass{article}

\usepackage{amsmath}
\usepackage{amsfonts}
\usepackage{amssymb}
\usepackage{color}  
\usepackage{amsthm}
\usepackage{graphicx}
\usepackage{caption}
\usepackage{comment}
\usepackage{subcaption}
\usepackage{algpseudocode}
\usepackage{algorithm}
\usepackage[utf8]{inputenc}
\usepackage{tabularx}
\usepackage{multirow}
\usepackage{hyperref}
\usepackage[english]{babel}
\usepackage{setspace}
\usepackage{natbib}
\bibpunct[, ]{(}{)}{,}{a}{}{,}%

\renewcommand{\baselinestretch}{1.5}

\title{Detour Dual Optimal Inequalities for Column Generation with Application to Routing and Location}

\author{Julian Yarkony \textsuperscript{\rm 1, \rm 3}, Naveed Haghani\textsuperscript{ \rm2}, Amelia Regan\textsuperscript{ \rm 3}   \\[2ex] 

\textsuperscript{\rm 1} Laminaar Optimization Research Group, La Jolla, CA\\ 

\textsuperscript{\rm 2}University of Maryland, College Park, MD\\

\textsuperscript{\rm 3}University of California, Irvine, CA\\

\\
\\
\\
\\

}

\date{March 2021}
\date{August 2021}

\usepackage{natbib}
\usepackage{graphicx}

\begin{document}

\maketitle

\begin{abstract}
    We consider the problem of accelerating column generation (CG) for logistics optimization problems using vehicle routing as an example.  Without loss of generality, we focus on the Capacitated Vehicle Routing Problem (CVRP) via the addition of a new class of dual optimal inequalities (DOI) that incorporate information about detours from the vehicle routes.  These inequalities extend the Smooth-DOI recently introduced in the literature for the solution of certain classes of set-covering problems by CG. The Detour-DOI introduced in this article permit low cost swap operations between items on a given active route with items near to other items on that route to estimate (and bound) the values of the dual variables. Smooth-DOI in contrast only permit low cost swap operations between nearby items. The use of Detour-DOI permits a faster convergence of CG without weakening the linear programming relaxation. We then argue that these DOI can also be conveniently applied to single source capacitated facility location problems. These problems have been shown to be equivalent to a broad class of logistics optimization problems that include, for example telecommunication network design and production planning. The importance of developing vastly more efficient column generation solvers cannot be overstated. 
    Detour-DOI, which permit large numbers of columns to be expressed with a finite set of variables, contributes to this important endeavor. 
\end{abstract}

\renewcommand{\baselinestretch}{2.0}
\section{Introduction}

Expanded linear programming (LP) relaxations of exponential size in that of the compact formulation (in terms of number of variables) provide tight (or near tight) dual bounds for many important classes of optimization problems \citep{lubbecke2005selected}, including applications in logistics and transportation \citep{costa2019}, as well as computer vision/machine learning \citep{yarkony2020data,FlexDOIArticle,HPlanarCC}. Given the intractable size of those models, their solutions are handled via the dynamic generation of variables using column generation (CG) \citep{barnprice}.

Unfortunately, CG often suffers from degeneracy and dual oscillation that hijack its efficient solution in many practical situations. Several techniques have been developed to mitigate this undesired behavior via the stabilization of the dual variables associated with the linear relaxations solved throughout the execution process \citep{du1999stabilized,lubbecke2005selected, oukil2007stabilized}. One powerful class of stabilization techniques is referred to as dual optimal inequalities (DOI) \citep{ben2006dual,FlexDOIArticle}, that enforce valid constraints on the dual space to limit the oscillation of the dual solutions. DOI may dramatically shrink the dual space over which CG must search accelerating optimization.

Our work builds off the work on Smooth-Dual Optimal Inequalities (S-DOI) \citep{haghani2020smooth}, which show that dual variables must change smoothly across space in large classes of problems embedded on metric spaces, and show their efficiency on handling the Capacitated Vehicle Routing Problem (CVRP) \citep{laporte1983branch,costa2019} and the Single-Source Capacitated Facility Location Problem \citep{diaz2002branch}. These problems have been shown to be equivalent to a broad class of logistics optimization problems that include, for example telecommunication network design and production planning \citep{yano1984equivalence, thizy1994facility}. 

This paper introduces a new class of DOI relevant for the CVRP referred to as Detour-Dual Optimal Inequalities (DT-DOI). These expand the S-DOI and permit a wider range of columns to be expressed using the same number of actual columns. DT-DOI permit the restricted master problem to add detours to routes where a detour services an individual item. The use of Detour-DOI permits a faster convergence of CG without weakening the linear programming relaxation.

We organize this document as follows. In Section \ref{Sec_lit_rev} we review the literature on dual stabilization in the context of CG generally and DOI specifically. In Section \ref{SecRevCVRPCG} we provide formal mathematical description of the CVRP \citep{Desrochers1992,costa2019,baldacci2011new,righini2009decremental}, and the solution via CG. In Section \ref{Sec_S_doi} we review the S-DOI \citep{haghani2020smooth}. In Section \ref{Sec_Detour_DOI} we introduce our DT-DOI and the associated optimization algorithm. In Section \ref{exper} we provide experimental validation of our approach. In Section \ref{Sec_SSCFLP} we demonstrate how the DT-DOI can be applied to the Single Source Capacitated Facility Location Problem. In Section \ref{conc} we conclude and discuss extensions. In Table \ref{myTabNotatoin} we provide the acronyms used in this document.
%
%
\begingroup
\setlength{\tabcolsep}{10pt} 
\renewcommand{\arraystretch}{1.5} 
\begin {table}[H]
\caption {\bf{Acronyms}} \label{tab:title} 
\begin{center}
\begin{tabular}{|c|c|} 

 \hline
 CVRP & Capacitated Vehicle Routing Problem \\
 CG & Column Generation \\ 
 SSCFLP & Single Source Capacitated Facility Location Problem \\
 RMP & Restricted Master Problem \\
 DOI & Dual Optimal Inequalities \\
 S-DOI & Smooth Dual Optimal Inequalities \\
 D-DOi & Deep Dual Optimal Inequalities \\
 F-DOI & Flexible Dual Optimal Inequalities \\
 DT-DOI & Detour Dual Optimal Inequalities \\
 \hline
\end{tabular}

\end{center}
\label{myTabNotatoin}
\end{table}
\endgroup

\section{Literature Review}
\label{Sec_lit_rev}

The literature on column generation (CG) techniques is vast. Here we review the most relevant material on general stabilization methods, trust region based methods and dual optimal inequalities (DOI) used to accelerate CG. CG suffers from slow convergence when the number of items in a column becomes large. This tends to produce intermediate dual solutions that are sparse and do not share properties with those of dual optimal solutions. We now consider some methods designed to circumvent this difficulty.

\subsection{General Stabilization Methods}
  
Due to inherent instability, many methods of stabilization have been proposed over the last two decades. Du Merle et al formalized the idea of stabilized column generation in their 1999 paper of that name \citep{du1999stabilized}. That paper proposed a 3-piecewise linear penalty function to stabilize the column generation procedure. Ben Amor and Desrosiers later proposed a 5-piecewise linear penalty function for improved stabilization \citep{amor2006proximal}. Shortly after, Oukil et al use the same framework to attack highly degenerate instances of multiple-depot vehicle scheduling problems \citep{oukil2007stabilized}. Ben Amor et al later proposed a general framework for stabilized CG  algorithms in which a stability center is chosen as an estimate of optimal solution of the dual formulation \citep{amor2009choice}. Gonzio et al proposed a primal-dual CG method in which the sub-optimal solutions of the restricted master problem are obtained using an interior point solver that was proposed in an earlier paper by the first author \citep{gondzio1995hopdm}. They examine their solution method relative to standard CG and analytic center cutting plane method proposed by Babonneau et al. \citep{babonneau2006solving, babonneau2007proximal}. They found that while standard column generation is efficient for small problem instances, that the primal-dual column generation method achieved the best solutions on larger problems \citep{Gondzio2013New}.

\subsection{Trust Region Based Methods}
Trust regions based methods discourage \citep{du1999stabilized} or prevent \citep{marsten1975boxstep} the next dual solution from leaving the area around the best dual solution generated thus far (best refers to the dual solution with greatest Lagrangian bound).  This is done since the current set of columns provides little information regarding the Lagrangian bound of dual solutions dis-similar to any previously generated dual solution.  For reference we now define the Lagrangian bound, which provides a lower bound on the optimal solution to the master program given any dual solution generated during CG.  The Lagrangian bound is the objective of the restricted master problem (at that iteration) plus the reduced cost of the lowest reduced column that is generated immediately subsequently during pricing multiplied by the maximum number of columns used (in any integer solution).  At termination of CG the Lagrangian bound equals the objective of the optimal solution to the master problem since no column has negative reduced cost.  

Smoothing based approaches \citep{Pessoa2018Automation} are a simple class of trust region approaches that achieve excellent results in practice. Smoothing based approaches only differ from standard CG in the selection of the dual variables terms which pricing is done on. Specifically they use a convex combination of:
\\
1) The dual solution of the current restricted master problem and \\
2) The dual solution generated thus far with greatest Lagrangian bound. 

\subsection{Dual Optimal Inequalities}
Ben Amor et al introduced the notion of dual optimal inequalities (DOI), which provide provable bounds on the optimal dual solution and ensure that at each step of CG the dual solution lies in the corresponding space \citep{ben2006dual}. These are further explained in \citep{Gschwind2016Dual}. In the primal form DOI often correspond to slack variables on primal constraints \citep{haghani2020smooth}.  
DOI are provably inactive at termination of CG but may or may not be active in intermediate steps.

The research presented in this document explicitly draws on the work on previous work on DOI and therefore uses the same terminology. Similar techniques to DOI are used by other researchers but with different language. For example, Miranda and Garrido explore inequalities on the dual in an inventory routing problem \citep{miranda2008valid}.

\subsection{Deep Dual Optimal Inequalities}
Deep Dual Optimal Inequalities (D-DOI) were introduced by \citep{ben2006dual}. A DOI is one that all dual-optimal solutions to the master problem satisfy. In contrast a D-DOI  does not necessarily satisfy all of them. A D-DOI is merely guaranteed satisfy at least one dual optimal solution \citep{ben2006dual}. Sets of D-DOI are valid as long as they preserve at least one dual optimal solution to the original master problem in the feasible dual space. These are are explored in detail in \citep{Gschwind2017Stabilized} in which they are applied to a temporal knapsack problem; and  \citep{gschwind2019stabilized} in which they are applied to commodity constrained split delivery problem; and by \citep{koza2020deep} in which they are applied to capacitated fixed charge network design problem. That research found significant reductions in solution times, iterations and the number of columns generated when the D-DOIs were compared to unstabalized column generation. 
\subsection{Flexible Dual Optimal Inequalities}
Flexible-Dual Optimal Inequalities (F-DOI) have been proposed recently, \citep{FlexDOIArticle,haghani2020relaxed,haghani2020smooth} and exploit the following observation:  the change in the cost of a column induced by removing a small number of items is often small, and can be easily bounded. Such bounds are column specific and exploit properties of the problem domain. For example in a set covering problem, in the primal form F-DOI correspond to additional variables that provide rewards for over-covering items. These rewards are set such that at optimality they are not used, but prior to termination of CG they have an important role. They can be understood as adding to the RMP all columns consisting of subsets of columns in the RMP. The costs of these columns provide upper bounds on the true cost of the column. F-DOI and their predecessors, varying and invariant DOI \citep{yarkony2020data}, provide considerable speed-ups and also make CG more robust to the specific selection of optimization parameters.   

\subsection{Family Column Generation}

Haghani et al recently introduced a new way to accelerate the convergence of column generation when applied to set-covering-based formulations by stabilizing dual optimization \citep{haghani2021family} called Family Column Generation (FCG). FCG can be applied on any set covering CG based formulation and accelerates optimization.  

FCG seeks to solve a restricted master problem over the set of all columns in the union of the families of columns in the restricted master problem. The so-called family of a column $l$ is a subset of the columns, which includes $l$, that is easy to price over, and is application specific. A specific solver is designed so that the set of all columns in each family need not be explicitly enumerated.   

\section{CVRP and the Column Generation Solution}
\label{SecRevCVRPCG}

We now consider the Capacitated Vehicle Routing Problem (CVRP) which is defined as follows. We are given a set of items, and a starting/ending depot embedded in a metric space. Each item is associated with an integer demand and the vehicles have a common capacity. We seek to partition the items into ordered lists of items called routes each serviced by a unique vehicle so as to minimize the total distance traveled while ensuring that no vehicle services more demand than it has capacity.  The number of vehicles used is bounded.  

We now consider the mathematical description of CVRP. We use $N$ to denote the set of items, which we index by $u$.  We define $N^+$ to be $N$ augmented with the starting depot, and ending depot which are denoted $-1,-2$ respectively (and are typically the same place). Items and depots are embedded on a metric space, and hence distances between them satisfy the triangle inequality. We use $c_{uv}$ to denote the distance between any pair of $u,v$ each of which lie in $N^+$.  
Each item is associated with an integer demand $d_u$, which lies in the set of unique demands $\mathcal{D}$ where $D_0$ is the capacity of a vehicle.  A route is feasible if it satisfies the following.  
\begin{itemize}
    \item 
The route starts and ends at the starting/ending depot respectively; 
\item The route visits an item no more than once.
\item The route services total demand that does not exceed $D_0$.
\end{itemize}
We denote the set of routes with $\Omega$, which we index by $l$. We describe the mapping of items to routes using  $a_{ul} \in \{0,1\}$ where $a_{ul}=1$ if and only if route $l$ contains item $u$ for any $u \in N$. For short hand we use  $N_{l} =\{ \forall u \in N \mbox{ s.t. } a_{ul}=1\}$ meaning that $N_l$ is the set of items serviced by  $l$. For any $u,v$ pair where each lie in $N^+$ we set $a_{uvl}=1$ IFF $v$ follows $u$ immediately in route $l$ and otherwise set $a_{uvl}=0$. The cost of a route is denoted $c_l$ is defined as the total distance traveled,  which we write formally as follows.
\begin{align}
    c_l=\sum_{\substack{u \in N^+\\ v \in N^+ }}c_{uv}a_{uvl}\quad  \forall l \in \Omega
\end{align}
The constraint that a route services no more demand than $D_0$ is written below using $d_u$ to denote the demand of item $u$.
\begin{align}
    \sum_{u \in N}d_u a_{ul}\leq D_0 \quad \forall l \in \Omega
\end{align}
A set of routes provides a feasible solution to CVRP if it services every item at least once and uses no more than $K$ routes where $K$ is the number of vehicles available.  We describe a solution to CVRP using decision variables $\theta_l \in \{0,1\}$ where $\theta_l=1$ indicates that route $l$ is selected and otherwise $\theta_l=0$.  We write the selection of the optimal solution below as follows with annotation describing the equations subsequently.  
\begin{subequations}
\label{origILP}
\begin{align}
    \min_{\theta \in \{0,1\} }\sum_{l \in \Omega}c_l\theta_l \label{ILPOBJ}\\
    \sum_{l \in \Omega}a_{ul}\theta_l \geq 1 \quad \forall u \in N \label{ILPCover}\\ 
    \sum_{l \in \Omega}\theta_l \leq K \label{ILPPack}
\end{align}
\end{subequations}
In \eqref{ILPOBJ} we seek to minimize the total cost of the routes selected.  In \eqref{ILPCover} we ensure that every item is covered at least once.  We should note that in any optimal solution that each item is covered exactly once since $c_l$  
increases as $N_l$ grows.  In \eqref{ILPPack} we enforce that no more than $K$ routes are used.  

We solve \eqref{origILP} via solving the linear programming (LP) relaxation, which is referred to as the master problem (MP). This LP relaxation is very tight in practice and can be tightened using branch-cut-price \citep{barnprice,lubbecke2005selected}. We write the MP below with dual variables written in brackets ($[ ]$) after the equations.  
\begin{subequations}
\label{primal_master}
\begin{align}
   \min_{\theta \geq 0 }\sum_{l \in \Omega}c_l\theta_l \\
    \sum_{l \in \Omega}a_{ul}\theta_l \geq 1 \quad  \forall u \in N  \quad [\pi_u]\\ 
    \sum_{l \in \Omega}\theta_l \leq K \quad [-\pi_0]
\end{align}
\end{subequations}
We write the dual form of \eqref{primal_master} below.  
\begin{subequations}
\label{dual_master}
\begin{align}
   \max_{\pi \geq 0}-K\pi_0+\sum_{u \in N}\pi_u \\
    c_l +\pi_0-\sum_{u \in N}a_{ul}\pi_u\geq 0\quad \quad \forall l \in \Omega 
\end{align}
\end{subequations}
Since the set of routes ($\Omega$) grows exponentially in the number of items we cannot trivially solve \eqref{primal_master}. Instead column generation (CG) is employed to solve \eqref{primal_master}. CG constructs a sufficient subset of $\Omega$ denoted $\Omega_R$ s.t. solving \eqref{primal_master} over $\Omega_R$ provides an optimal solution to \eqref{primal_master} over $\Omega$. To construct $\Omega_R$, we iterate between \textbf{(1)} solving \eqref{primal_master} over $\Omega_R$, which is referred to as the restricted master problem (RMP) and \textbf{(2)} identifying variables with negative reduced cost, which are then added to $\Omega_R$. Typically the lowest reduced cost column is generated.  We write the selection of this column as optimization below using $\bar{c}_l$ to denote the reduced cost of column $l$.  
\begin{subequations}
\label{pricing}
\begin{align}
    \min_{l \in \Omega} \bar{c}_l \\
    \bar{c_l}=c_l+\pi_0-\sum_{u \in N}a_{ul}\pi_u \label{redCostForm}
\end{align}
\end{subequations}
The operation in \eqref{pricing} is referred to as pricing.  One or more negative reduced cost columns are generated during pricing. CG terminates when pricing finds no column with negative reduced cost.  This certifies that CG has produced the optimal solution to \eqref{primal_master}.  We write pricing as an integer linear program (ILP) in Appendix \ref{pricingILP} though it is often solved with a resource constrained shortest path solver \citep{costa2019}.  We terminate CG when no $l$ in $\Omega$ has negative reduced cost.
CG initializes $\Omega_R$ with a heuristically generated feasible integer solution or using artificial variables that have prohibitively high cost but ensure a feasible solution.  
In Alg \ref{basicCG} we describe CG in pseudo-code.
\begin{algorithm}[!b]
 \caption{Basic Column Generation}
\begin{algorithmic}[1] 
\State $\Omega_R\leftarrow $ from user
\label{Zline_rec_input_start}
\Repeat
\label{Zline_outer_start}%
\State  Solve for $\theta,\pi$  using \eqref{primal_master},\eqref{dual_master} over $\Omega_R$
\label{Zline_solve_FRMP}
\State $l_* \leftarrow \min_{l \in \Omega}\bar{c}_l$
\State $\Omega_R \leftarrow \Omega_R \cup l_*$
 \Until{$\bar{c}_{l_*} \geq 0$}
 \State Return last $\theta$  generated.  \label{ZreturnSol}
\end{algorithmic}
\label{basicCG}
\end{algorithm}

\section{Smooth Dual Optimal Inequalities}
\label{Sec_S_doi}

Smooth-Dual Optimal Inequalities (S-DOI) \citep{haghani2020smooth} exploit the fact that problems in operations research are often defined on metric spaces to provide bounds on dual variables. S-DOI enforce that dual variables have to change smoothly across space. S-DOI can be applied to CVRP and of other problems including the Single Source Capacitated Facility Location Problem.  S-DOI rely on a proof that an optimal dual solution obeys the following (using $\rho_{uv}$ defined subsequently).  
\begin{align}
\label{S_DOI_eq}
    \rho_{uv}\geq \pi_v -\pi_u \quad  \forall u \in N,v \in N \quad \mbox{s.t.} \quad  d_u\geq d_v 
\end{align}
Here $\rho_{uv}$ is an upper bound on the increase of the cost of a column containing $u$ but not $v$ when $u$ get replaced by $v$.  In the context of CVRP $\rho_{uv}=2c_{uv}$ where $c_{uv}$ is the distance from $u$ to $v$.  Here $2c_{uv}$ corresponds to making a detour at $u$ traveling to $v$ then traveling back to $u$ then continuing the route.  In this context $u$ is not serviced but $v$ is serviced. 

In CG optimization using S-DOI, \eqref{S_DOI_eq} is enforced at each step of CG optimization.  The application of S-DOI in \citep{haghani2020smooth} produces significant improvements when compared to standard CG in the Single Source Capacitated Facility Location Problem \citep{diaz2002branch,holmberg1999exact}. 

The primal form of S-DOI is written below using additional terms $\omega_{uv}$ to denote the decision variables for swapping item $u$ for $v$.  We use $S$ to denote the set of $u \in N,v\in N$ for which $d_u\geq d_v$, and hence $u$ can be swapped for $v$ while preserving the feasibility any given column containing $u$.
\begin{subequations}
\label{SDOIOPT}
\begin{align}
\min_{\substack{\theta \geq 0\\ \omega \geq 0} }\sum_{l \in \Omega}c_l\theta_l+\sum_{uv \in S}\rho_{uv}\omega_{uv} \\
    \sum_{\substack{v \in N\\ uv \in S}}-\omega_{uv}+\sum_{\substack{v \in N \\ vu \in S}}\omega_{vu}+\sum_{l \in \Omega}a_{ul}\theta_l \geq 1 \quad \forall u \in N  \\ 
    \sum_{l \in \Omega}\theta_l \leq K
\end{align}
\end{subequations}
The dual form of \eqref{SDOIOPT} is identical to \eqref{dual_master} except that
\eqref{S_DOI_eq} is enforced. CG proceeds as in Alg \ref{basicCG} but where \eqref{SDOIOPT} and its dual are solved to generate a $\theta,\pi$ pair. Pricing is unmodified for CG optimization. Typically $\rho$ is offset by a tiny positive constant so as to guarantee that at termination of CG that no $\omega$ terms are active in an optimal primal solution at termination of CG.  This offset is applied in the context of other DOI including Flexible-DOI \citep{FlexDOIArticle}.  
\section{DT-DOI:  Detour Dual Optimal Inequalities}
\label{Sec_Detour_DOI}
In this section we introduce the Detour-DOI (DT-DOI), which build off of the S-DOI, and further accelerate optimization without loosening the relaxation. S-DOI provide for the ability to swap nearby items at low cost. However we would like to swap distant items at low cost, which is not permitted by S-DOI. This comes up in the following context. Consider that we have a route in the RMP (in $\Omega_R$) servicing five items in Southern California and six items in  Montreal Canada; all items are of equal demand of one; and the capacity bound is 11 units (meaning $D_0=11$). There are a total of eight units each in Southern California and Montreal Canada. The route proceeds as follows:  (1) The vehicle starts at the starting depot; (2) then travels Southern California where it services its associated items there; (3) then  travels to Montreal where it services its items there; (4) then travels to the ending depot. The use of S-DOI permits us at low cost to alter this route to service any five items in Southern California AND any six items in Montreal because proximate items have nearly identical service cost. However we are unable to express the desire to have this route service a different number of items in each location at low cost. For example we can not express that this route should service eight items in Southern California and three in Montreal without using $\omega_{uv}$ terms associated with large $\rho_{uv}$ that correspond to swapping an item in Montreal with one in Southern California. We seek to correct this issue with the DT-DOI.  

We now describe an alternative primal master formulation for the CVRP. This formulation is augmented with additional primal variables/constraints.  
In Appendix \ref{sec_dualProof} we prove that the relaxation used in DT-DOI is exactly as tight as the standard set cover relaxation in \eqref{primal_master} and hence at termination of CG (at optimality) results produced are identical to standard CG. However prior to convergence of CG this relaxation produces a value between those of the MP and the standard RMP and thus permits CG to converge faster.  

We now provide additional notation used to express DT-DOI.  
We use decision variable $y_{ul} \in [0,1]$ to indicate the decision to make detour to service item $u$ on route $l$. This detour is made at the closest item (or depot) on the route $l$ to $u$. If the detour is made at item $v$ then the corresponding vehicle reaches $v$ then visits $u$ then returns to $v$ then proceeds on its route.  We set $y_{ul}=1$ if we make a detour on route $l$ to service item $u$. Note that $y_{ul}$ may be set to one or zero independent of $a_{ul}$.  If $y_{ul}=a_{ul}=1$ then the item $u$ is serviced as it normally would be on route $l$. We can understand this as making a detour of distance zero. Similarly if $y_{ul}=0$ while $a_{ul}=1$ then the route visits $u$ without servicing $u$.  The cost of this detour is written as $c_{ul}$ where $c_{ul}=2*$ the minimum distance from an item $v \in N_l$ (or the starting/ending depot) to $u$. We write this formally below (recall -1,-2 denote the start/end depot respectively; which are typically the same place).
\begin{align}
\label{DDOICVRPdef}
    c_{ul}=2\min_{v \in N_l+(-1)+(-2)} c_{uv}
\end{align}
We use decision variables $\psi_l \in [0,1]$ to indicate the selection of a route using detours.  We set $\psi_l=1$ to select route $l$ augmented with the ability to make detours.  We use $\mathcal{D}$ to denote the set of unique demands of items.  We define $D_{dl}$ to be the number of items in route $l$ of demand at least $d$ meaning $D_{dl}=\sum_{u \in N}[d_u\geq d]a_{ul} $ for all $l\in \Omega,d \in \mathcal{D}$.  We now write our primal master problem below with exposition afterwords.%
%
%
\begin{subequations}
\label{DT_DOI_master_primal}
\begin{align}
    \min_{\substack{\theta \geq 0 \\ \psi \geq 0 \\ y \geq 0}}\sum_{l \in \Omega}c_l\theta_l+\sum_{l \in \Omega}c_l\psi_l+\sum_{\substack{u \in N\\ l \in \Omega}}c_{ul}y_{ul} \label{DDOIOBJ}\\
    \sum_{\substack{ l \in \Omega}}y_{ul}+\sum_{l \in \Omega}a_{ul}\theta_l \geq 1 \quad \quad \forall u \in N \quad [\pi_u] \label{DDOIcover}\\
    y_{ul}\leq \psi_{l} \quad \forall u \in N, l \in \Omega \quad [-\pi_{ul}]\label{l_con_1}\\
    \sum_{u \in N}[d_u\geq d]y_{ul} \leq D_{dl}\psi_l\quad \forall l \in \Omega,d \in \mathcal{D} \quad [-\pi_{dl}]\label{l_con_2}\\
    \sum_{l \in \Omega}\theta_l+\sum_{l \in \Omega}\psi_l \leq K \quad [-\pi_0]\label{capDDOI}
\end{align}
\end{subequations}

In \eqref{DDOIOBJ} we minimize the total cost of the routes and detours taken. In \eqref{DDOIcover} we enforce that an item is covered at least once by a detour or a route without detours. In \eqref{l_con_1} we enforce that we only take a detour on a given route if that route is selected.  In \eqref{l_con_2} we enforce that each route is feasible with respect to demand. Specifically we enforce that a route does not service more items with demand greater than or equal to any given amount than it services in the original route (where the original route is associated with $\theta_l$). In \eqref{capDDOI} we enforce that no more than $K$ vehicles are used.  We write the dual form of \eqref{DT_DOI_master_primal} below.  
\begin{subequations}
\label{DT_DOI_master_dual}
\begin{align}
    \max_{\pi \geq 0}-K\pi_0+\sum_{u \in N}\pi_u\\
 c_l+\pi_0-\sum_{u \in N}a_{ul}\pi_u\geq 0\quad \forall l \in \Omega \label{dual_con_1}\\
    c_l-\sum_{d \in \mathcal{D}}D_{dl}\pi_{dl}+\pi_0-\sum_{\substack{u \in N\\ l \in \Omega}}\pi_{ul}\geq 0 \quad \forall l \in \Omega  \label{eqdualTricky}\\
    c_{ul}+\pi_{ul}-\pi_u+\sum_{d\in \mathcal{D}}[d_u\geq d]\pi_{dl}\geq 0\quad \quad \forall l \in \Omega, u \in N\label{dual_con_3}
\end{align}
\end{subequations}
We use CG to solve \eqref{DT_DOI_master_primal}. In this case we replace $\Omega$ with $\Omega_R$. Observe that pricing over $\Omega$ is unaffected and is still done using \eqref{pricing}. \textbf{Note that no pricing is done over $y_{ul}$,or $\psi_l$.}  This is because once $l$ is added to $\Omega_R$ via \eqref{pricing} we add $\psi_l$ and $y_{ul} \; (\forall u \in N)$ to consideration in the RMP.

If it becomes the case that the RMP becomes difficult to solve we can decrease the number of $\psi_l$ terms, and or $y_{ul}$ terms considered.  Specifically we can choose to use only $\psi_{l},y_{ul}$ terms that we think would improve the objective. We did not consider this in experiments since the solution to the RMP was such a small part of computation time in our experiments relative to pricing time.   We consider an equivalent version of \eqref{DT_DOI_master_primal} with a reduced number of variables in Appendix \ref{APP_removeThetaAppend}. 
\section{Experimental Analysis}
\label{exper}
In this section we consider experimental evidence supporting the use of  Detour-DOI (DT-DOI) for accelerating column generation (CG) on the Capacitated Vehicle Routing Problem (CVRP).  We compared the performance of our DT-DOI against Smooth-DOI (S-DOI) and un-stabilized CG.  We did comparisons with regards to time and iterations to solve the linear programming relaxation of set cover for the CVRP.  

To provide fair comparisons we provided a ``vanilla" implementation of pricing.  We generate one column at each iteration of CG; where this column is the lowest reduced cost column.  Hence no heuristics are used during pricing.  To do this we solve the integer linear programming (ILP) formulation of pricing described in Appendix \ref{pricingILP}.  

We considered fourteen instances generated randomly instances of the following form.  Each problem instance is associated with $40$ items of demand one and $5$ vehicles of capacity $10$ each.  Each item and the starting depot is assigned a random integer position on grid of size 100 by 100.  Distances between items (and depots) are computed based on the L2 distance rounding up to the nearest integer.  The starting depot is located at the same place as the ending depot.

To solve the restricted master problem during the course of CG optimization we used the basic MATLAB linear programming solver with default options. For pricing we used the CPLEX mixed integer linear programming solver with default options.  

We initialize CG with one artificial column for each item. This column uses no vehicle; has prohibitively high cost; and covers the corresponding item.  

\begin{figure}[!htpb]
\centering
	\includegraphics[width=0.49\linewidth]{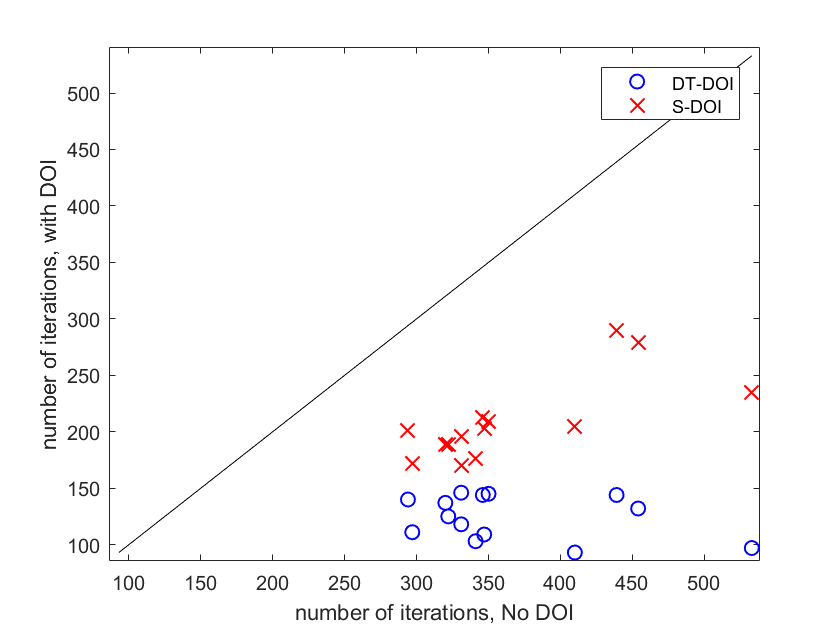}
	\includegraphics[width=0.49\linewidth]{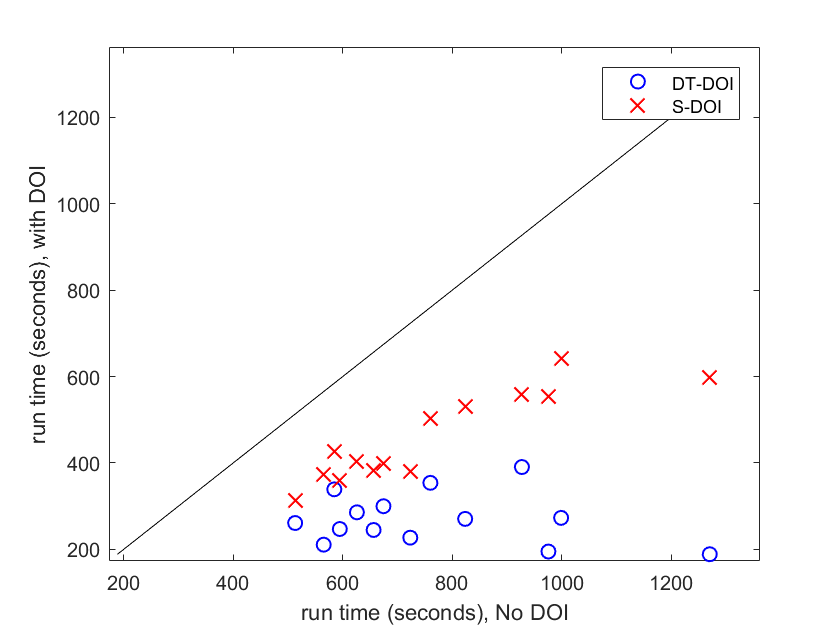}
	\caption[Detour-DOI results]{
		\textbf{(Left):  }Comparative iterations required between using stabilization and using no stabilization for all problem instances.  Each blue dot describes the performance of DT-DOI on a single problem instance.  Its x-coordinate is the number of iterations DT-DOI with CG took to solve that problem instance; while the y-coordinate is the number of iterations required to solve that problem instance with un-stabilized CG.  Red dots compare CG with S-DOI to un-stabilized CG.  The black line plots the line y=x so as to provide a baseline for improvement.  The further below the line the greater the improvement achieved.  
		\textbf{(Right):} Comparative run times between using stabilization and using no stabilization for all problem instances.  Lines and dots describe the same terms as (Left) except considering run time.}
	\label{fig:CVRP_results}
\end{figure}

We provide convergence time and iteration count aggregate plots in Tables \ref{table:CVRP_results_time} and Table \ref{table:CVRP_results_iter} respectively.  These show how long CG with DT-DOI, CG with S-DOI, and un-stabilized CG take to solve each problem instance. This data is aggregated in plots in Fig \ref{fig:CVRP_results}.

\begin{figure}[!hbtp]
	\includegraphics[width=0.49\linewidth]{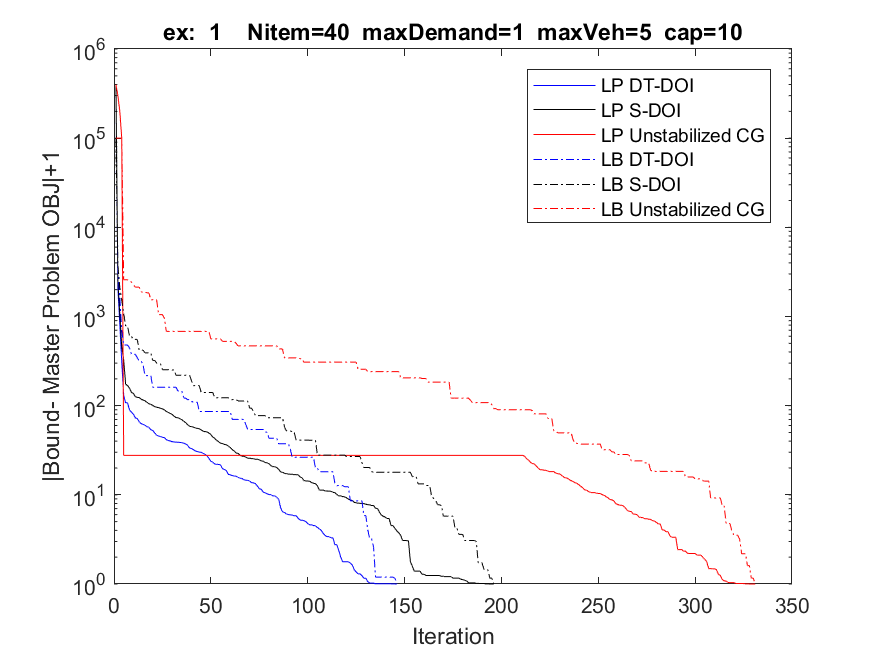}
	\includegraphics[width=0.49\linewidth]{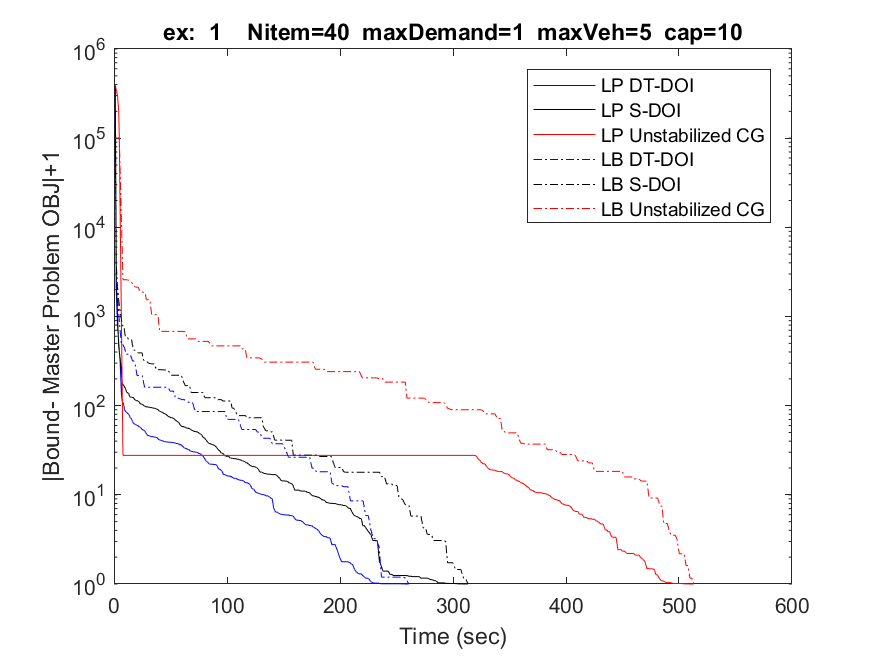}\\
	\includegraphics[width=0.49\linewidth]{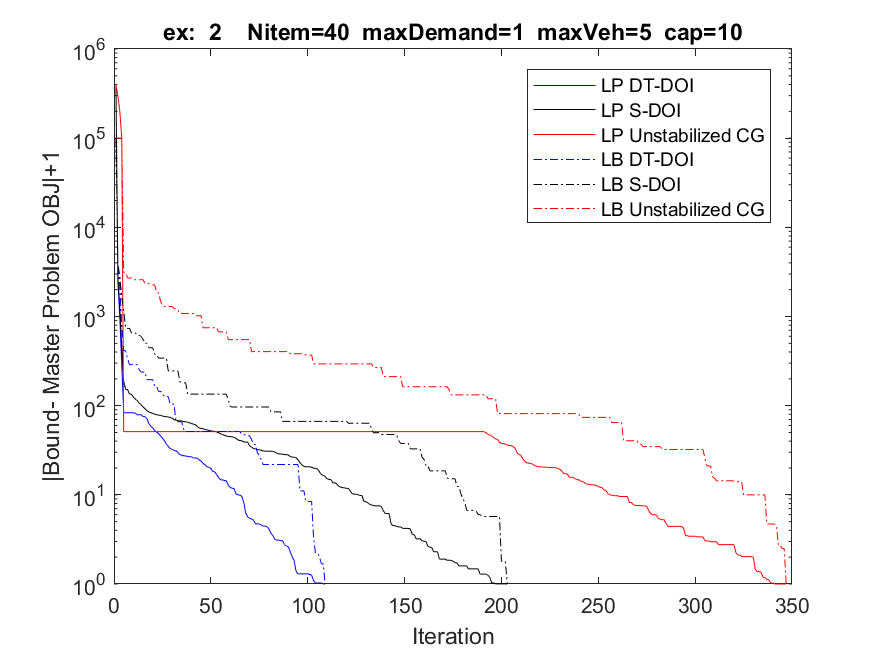}
	\includegraphics[width=0.49\linewidth]{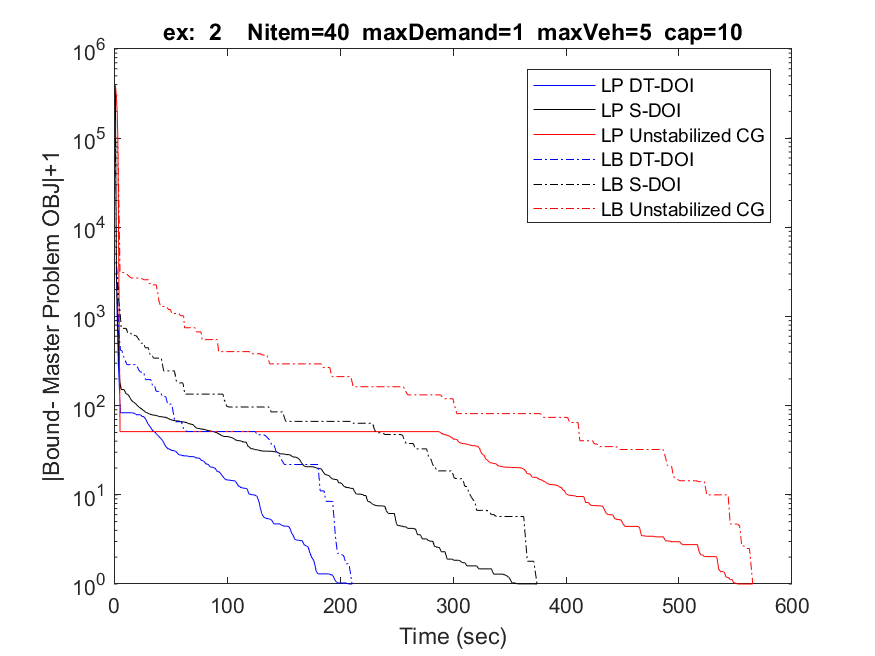}\\
	\includegraphics[width=0.49\linewidth]{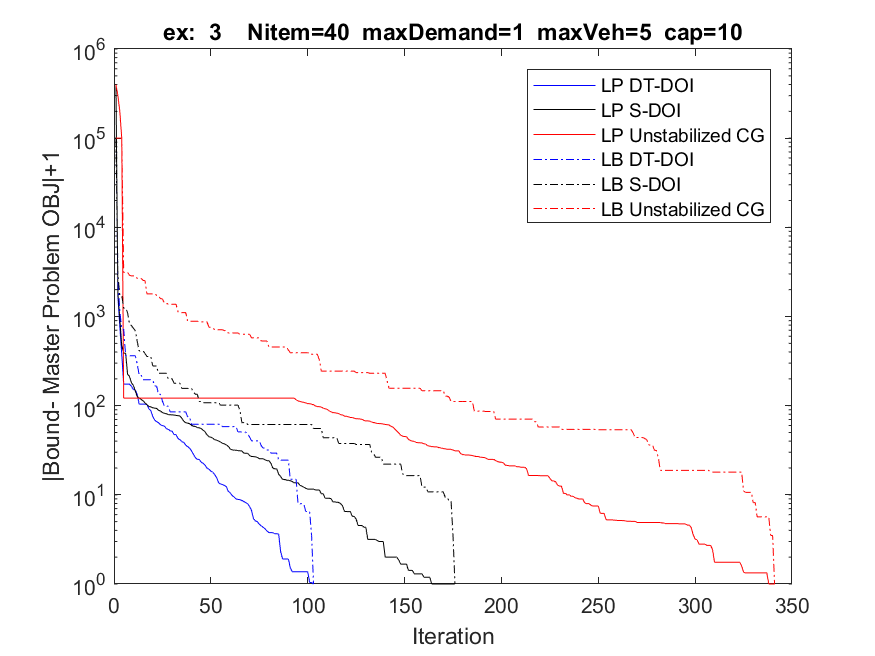}
	\includegraphics[width=0.49\linewidth]{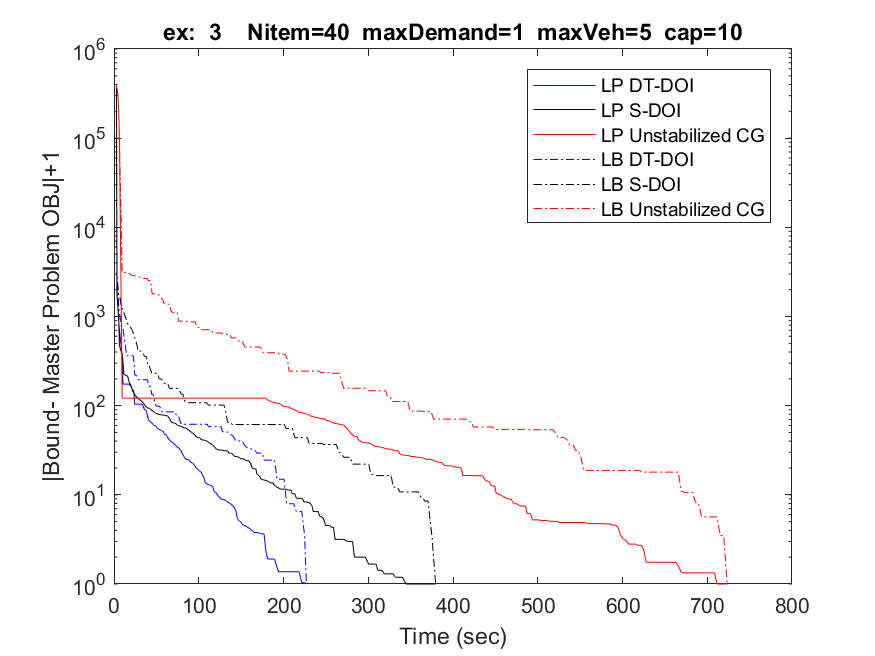}\\
	\caption{Results as a function of iteration/time.  The left side provides the results with respect to iteration and the right side for time on the same instance. We add one to the difference of all bounds and MP values  which allows us to use the semilog scale.  
	}
	\label{fig:iter1}
\end{figure}

\begin{figure}[!hbtp]
	\includegraphics[width=0.49\linewidth]{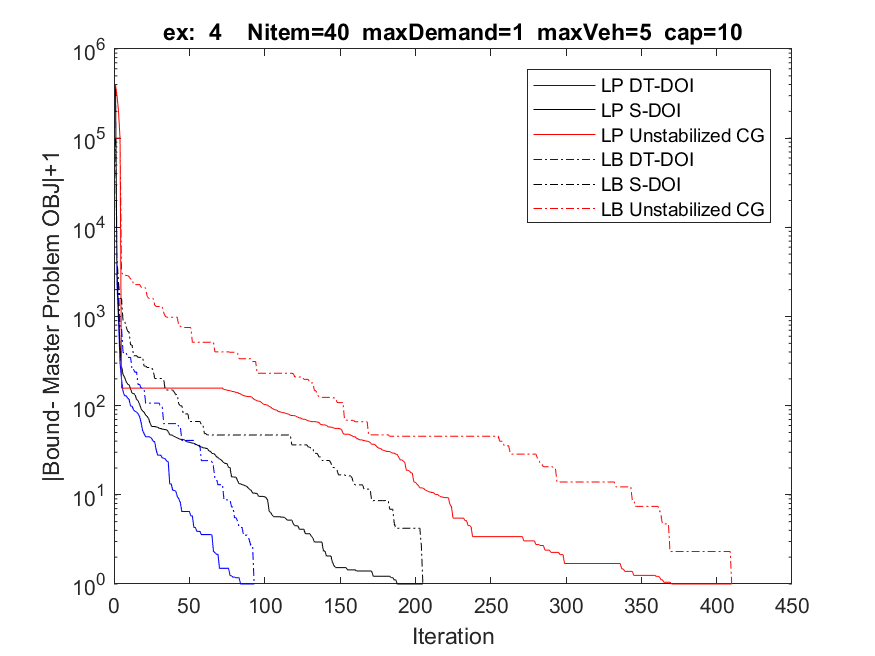}
	\includegraphics[width=0.49\linewidth]{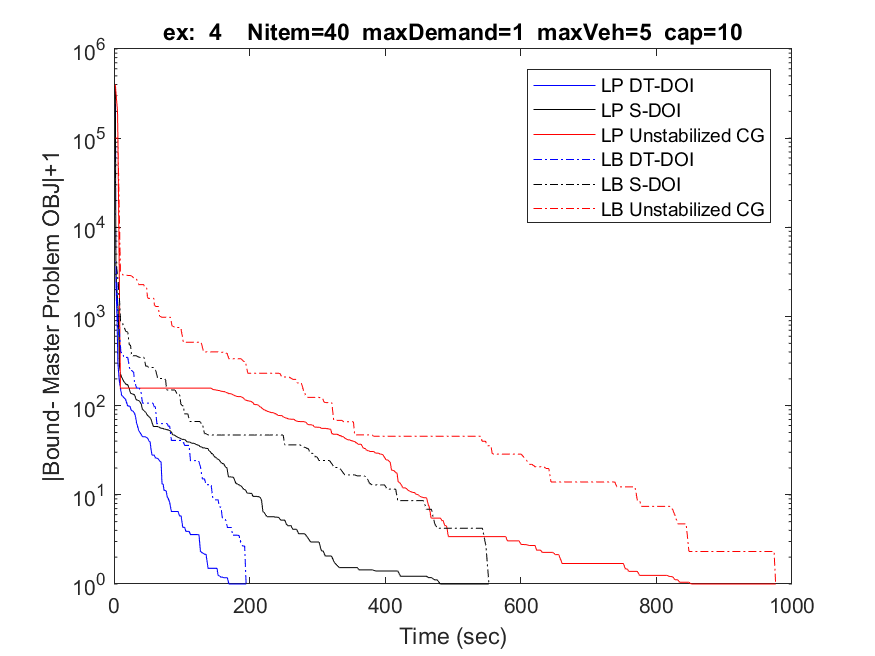}\\
	\includegraphics[width=0.49\linewidth]{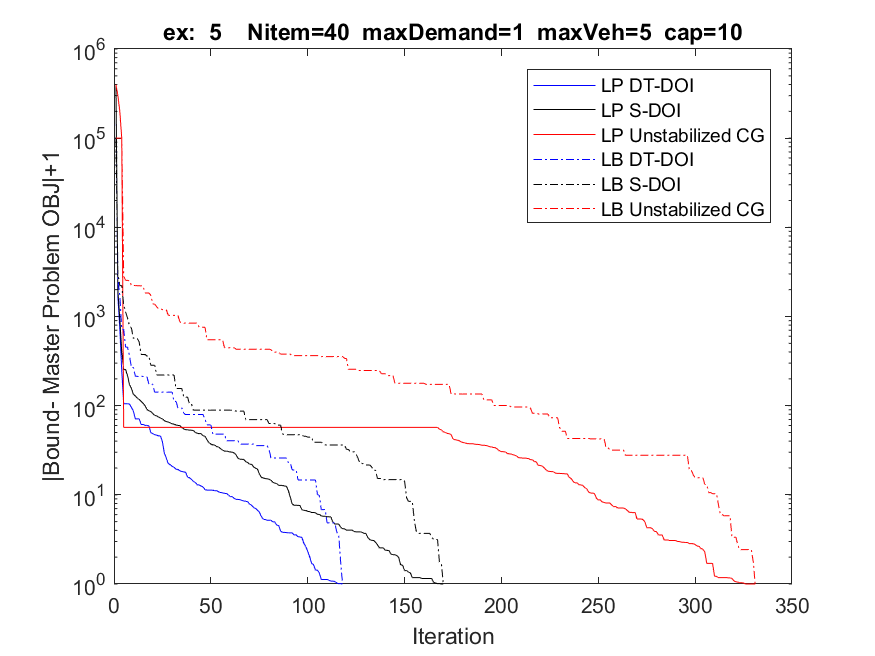}
	\includegraphics[width=0.49\linewidth]{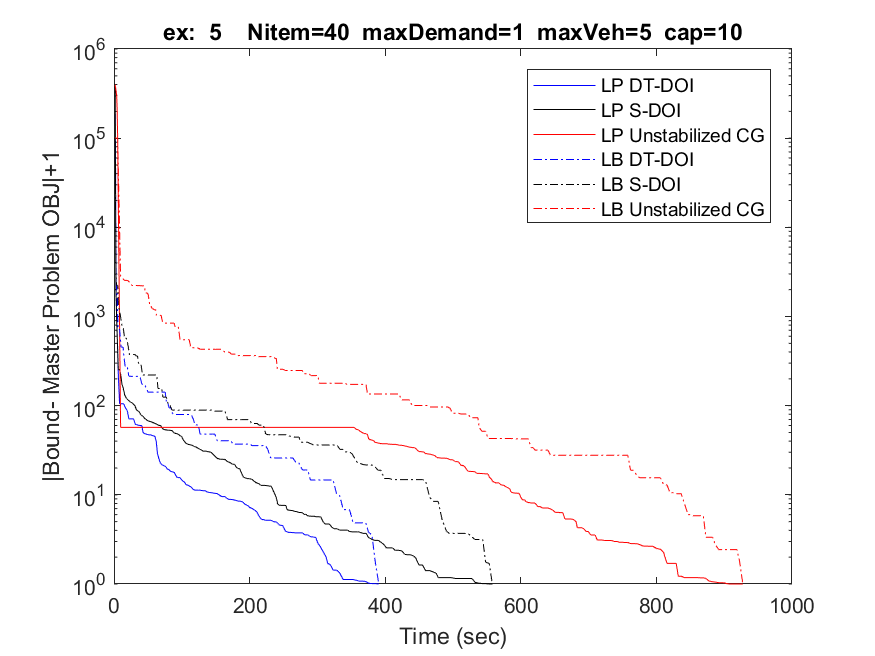}\\
	\includegraphics[width=0.49\linewidth]{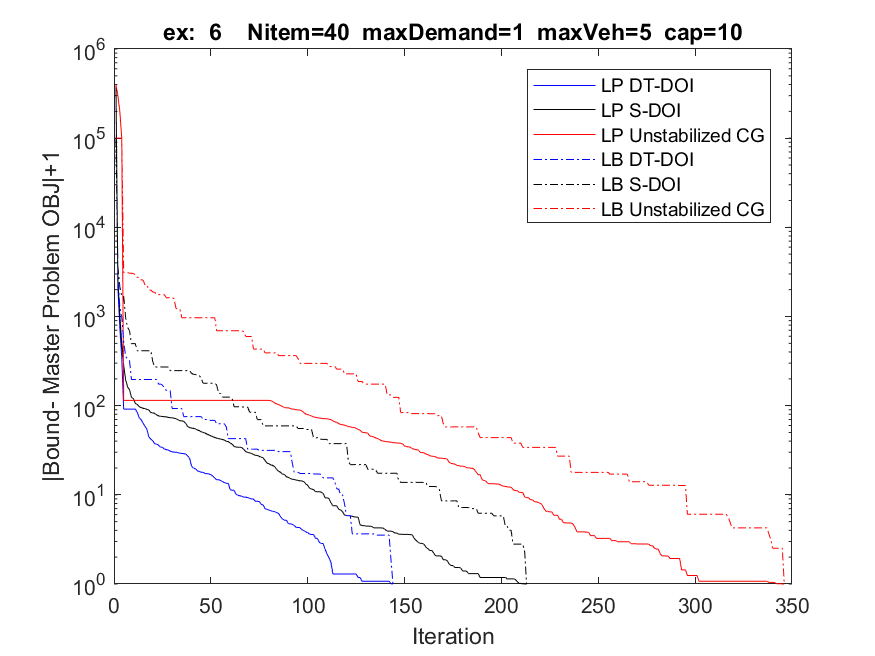}
	\includegraphics[width=0.49\linewidth]{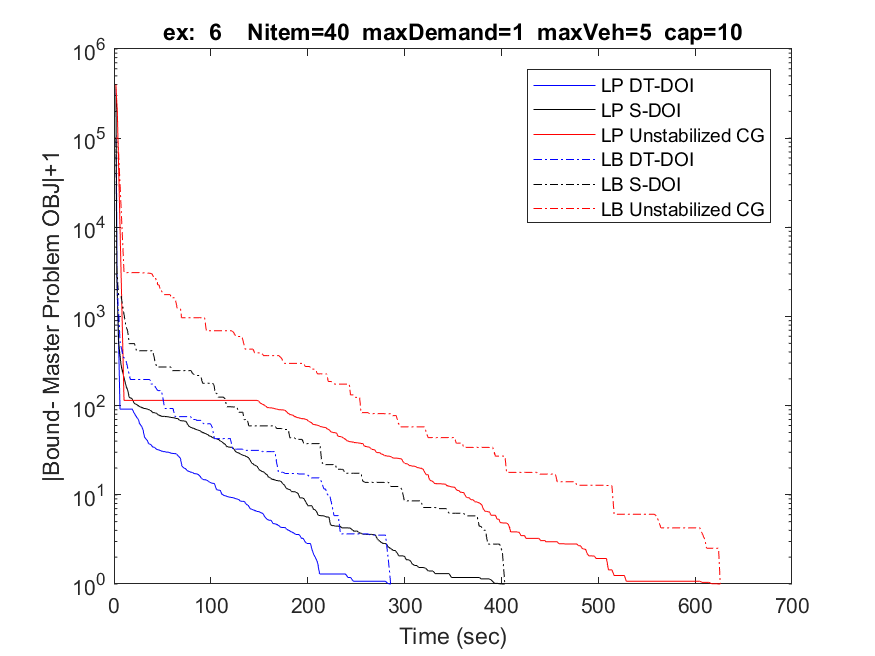}\\
	\caption{Results as a function of iteration/time.  The left side provides the results with respect to iteration and the right side for time on the same instance.  We add one to the difference of all bounds and MP values  which allows us to use the semilog scale.  
	}
	\label{fig:iter2}
\end{figure}

We provide results on our problem instances in Figs \ref{fig:iter1} and \ref{fig:iter2} in this section. Because the results are very similar, we have moved additional results in Figs \ref{fig:iter3}, \ref{fig:iter4} and \ref{fig:iter5} into the Appendix \ref{addFig}. Each figure describes performance on a separate problem instance with regards to iterations or time. We display the value of the linear programming relaxation and the corresponding Lagrangian bound as a function of time/iteration. The Lagrangian bound provides a lower bound on the optimal LP relaxation \citep{lubbecke2005selected}. 

The Lagrangian bound is the value of the linear program plus the value of the lowest reduced cost column times the number of vehicles $K$. For each algorithm we provide the value of the greatest Lagrangian bound computed thus far when displaying it, which we refer to as LB for lower bound.  
For both the LP values and the LB we display the difference between their value and the optimal value of the master problem (MP) as an absolute value (plus one). The plus one allows us to use the semilog scale.  

Empirically we observe large speed ups with regards to both timing and iteration for both S-DOI and DT-DOI over un-stabilized CG.  We furthermore see large speed ups for DT-DOI over S-DOI.  Computation time is dominated by pricing for all algorithms.

\begin{table}[tpb]
	\centering
	\scalebox{0.85}{
	\begin{tabular}{|c|c c c|c c|} 
		\hline
		& \multicolumn{3}{c|}{\bf Time} & \multicolumn{2}{c|}{\bf Speedup}\\
		\bf Instance & \bf un-stabilized CG & \bf DT-DOI & \bf S-DOI & \bf DT-DOI & \bf S-DOI \\
		\hline
        1 & 513.2 & 260.9 & 313.2 & 2.0 & 1.6 \\
        2 & 565.3 & 210.7 & 374.4 & 2.7 & 1.5 \\
        3 & 723.5 & 227.0 & 379.5 & 3.2 & 1.9 \\
        4 & 975.8 & 194.9 & 552.9 & 5.0 & 1.8 \\
        5 & 927.4 & 390.6 & 558.0 & 2.4 & 1.7 \\
        6 & 626.0 & 285.7 & 403.5 & 2.2 & 1.6 \\
        7 & 594.7 & 247.0 & 360.2 & 2.4 & 1.7 \\
        8 & 656.4 & 244.9 & 382.1 & 2.7 & 1.7 \\
        9 & 674.5 & 299.8 & 397.6 & 2.3 & 1.7 \\
        10 & 823.8 & 270.5 & 529.8 & 3.0 & 1.6 \\
        11 & 999.0 & 272.9 & 641.9 & 3.7 & 1.6 \\
        12 & 760.2 & 353.9 & 503.9 & 2.1 & 1.5 \\
        13 & 584.7 & 338.9 & 425.5 & 1.7 & 1.4 \\
        14 & 1270.6 & 188.6 & 598.0 & 6.7 & 2.1 \\
        \hline
        mean & 763.9 & 270.4 & 458.6 & 3.0 & 1.7 \\
        median & 699.0 & 265.7 & 414.5 & 2.5 & 1.6 \\
		\hline
	\end{tabular}}
	\caption{CVRP runtime results}
	\label{table:CVRP_results_time}
\end{table}

\begin{table}[tpb]
	\centering
	\scalebox{0.85}{
	\begin{tabular}{|c|c c c|c c|} 
		\hline
		& \multicolumn{3}{c|}{\bf Iterations} & \multicolumn{2}{c|}{\bf Iteration Speedup}\\
		\bf Instance & \bf un-stabilized CG & \bf DT-DOI & \bf S-DOI & \bf DT-DOI & \bf S-DOI \\
		\hline
        1 & 331 & 146 & 196 & 2.3 & 1.7 \\
        2 & 347 & 109 & 203 & 3.2 & 1.7 \\
        3 & 341 & 103 & 176 & 3.3 & 1.9 \\
        4 & 410 & 93 & 205 & 4.4 & 2.0 \\
        5 & 331 & 118 & 170 & 2.8 & 1.9 \\
        6 & 346 & 144 & 213 & 2.4 & 1.6 \\
        7 & 322 & 125 & 189 & 2.6 & 1.7 \\
        8 & 297 & 111 & 172 & 2.7 & 1.7 \\
        9 & 320 & 137 & 189 & 2.3 & 1.7 \\
        10 & 439 & 144 & 290 & 3.0 & 1.5 \\
        11 & 454 & 132 & 279 & 3.4 & 1.6 \\
        12 & 350 & 145 & 209 & 2.4 & 1.7 \\
        13 & 294 & 140 & 201 & 2.1 & 1.5 \\
        14 & 533 & 97 & 235 & 5.5 & 2.3 \\
        \hline
        mean & 365.4 & 124.6 & 209.1 & 3.0 & 1.8 \\
        median & 343.5 & 128.5 & 202 & 2.7 & 1.7 \\
		\hline
	\end{tabular}}
	\caption{CVRP iteration results}
	\label{table:CVRP_results_iter}
\end{table}


\clearpage
\section{Application to Single Source Capacitated Facility Location}
\label{Sec_SSCFLP}
In this section we apply our Detour-Dual Optimal Inequalities (DT-DOI) to the Single Source Capacitated Facility Location Problem (SSCFLP) \citep{diaz2002branch}. In SSCFLP we open a set of facilities with fixed capacity each of which services all demand of a subset of the customers  s.t. each customer is serviced. Note that the customers in the SSCFLP map to items in the CVRP. We seek to minimize the total cost of opening facilities and serving customers.  

In this problem we are provided with a set of customers $N$, which we index by $u$, and a set of facilities $\mathcal{F}$ indexed by $f$. Each facility has an integer capacity $D_f$, and each customer has an integer demand $d_u$, which lies in the set of demands $\mathcal{D}$.  Each facility has a fixed cost $c_f$ to open and a cost to service a customer $u$ denoted $c_{fu}$.  A solution is feasible if every customer is covered at least once and no facility is opened more than once.  The set of valid assignments (where an assignment is the analog of a route in CVRP) is denoted $\Omega$, and indexed by $l$. We use $a_{fl}=1$ if assignment $f$ is associated with $l$ and otherwise set $a_{fl}=0$.  Each assignment is associated with exactly one facility meaning $\sum_{f \in \mathcal{F}}a_{fl}=1$ for all $l\in \Omega$.  We use $a_{ul}=1$ to indicate that assignment $l$ services customer $u$, and otherwise set $a_{ul}=0$. The cost of a column is written $c_l$, which is defined as the fixed cost of opening facility $f$ plus the cost of servicing its respective customers. We formally define $c_l$ below.
\begin{align}
    c_l=\sum_{f \in \mathcal{F}}a_{fl}c_f+\sum_{\substack{f \in \mathcal{F}\\u \in N}}a_{fl}a_{ul}c_{fu} \quad \forall l \in \Omega
\end{align}
An assignment must satisfy all of the demand of its associated customers.  We write this constraint formally as follows.  
\begin{align}
   \sum_{u \in N}a_{ul}d_{u}\leq \sum_{f \in \mathcal{F}}a_{fl}D_f \quad \forall l \in \Omega
\end{align}
We use decision variable $\theta_l=1$ to indicate that assignment $l$ is selected.  We write the column generation (CG) master problem below with dual variables in ([]) and annotated the equations subsequently. 
\begin{subequations}
\label{SSCFLPmaster}
\begin{align}
    \min_{\theta \geq 0}\sum_{l \in \Omega} c_l\theta_l \label{sscflpOBJ}\\
    \sum_{l \in \Omega}a_{ul}\theta_l\geq 1 \quad \forall u \in N \quad [\pi_u]\label{fac_item_cover}\\
    \sum_{l \in \Omega}a_{fl}\theta_l\leq 1 \quad \forall f \in \mathcal{F}\quad [-\pi_f]  \label{fac_pack}
\end{align}
\end{subequations}
In \eqref{sscflpOBJ} we seek to minimize the total cost of the assignments selected. In \eqref{fac_item_cover} we enforce that each item is covered at least once.   In \eqref{fac_pack} we enforce that each facility is opened at most once.  

The corresponding pricing problem is solved independently for each $f \in \mathcal{F}$ as a knapsack problem.  This is written as an integer linear program (ILP) below using $x_u=1$ to select item $u$. 
\begin{subequations}
\label{SSCFLP_pricing}
\begin{align}
    \min_{x \in \{0,1\}} c_f+\pi_f+\sum_{u \in N}(c_{fu}-\pi_u)x_{u}\\
    \sum_{u \in N}d_ux_u\leq D_f
\end{align}
\end{subequations}

We now apply our Detour-DOI approach to SSCFLP.  One key difference is that cost of making a detour is dependent not on the customer in the assignment but only the facility used.  Thus  $c_{ul}=\sum_{f \in \mathcal{F}}a_{fl}c_{fu} \quad \forall l \in \Omega, u \in N$ and the objective term associated with $\psi_l$ can remove the costs for servicing customers in $N_l$ (which is the set of customers that $l$ services).  We refer to this cost as $\hat{c}_l$ which is defined as follows and includes only the $c_f$ term for the facility associated with $l$:  
$\hat{c}_l=\sum_{f \in \mathcal{F}}a_{fl}c_{f}$.
We now write the Detour-DOI formulation below using $y_{ul},D_{dl}$ as in \eqref{DT_DOI_master_primal}.\\

\begin{subequations}
\label{DT_DOI_master_primal3}
\begin{align}
    \min_{\substack{\theta \geq 0 \\ \psi \geq 0 \\ y \geq 0}}\sum_{l \in \Omega}c_l\theta_l+\sum_{l \in \Omega}\hat{c}_l\psi_l+\sum_{\substack{u \in N\\ l \in \Omega}}c_{ul}y_{ul} \label{DDOIOBJ3}\\
    \sum_{\substack{ l \in \Omega}}y_{ul}+\sum_{l \in \Omega}a_{ul}\theta_l \geq 1 \quad \quad \forall u \in N \quad [\pi_u] \label{DDOIcover3}\\
    y_{ul}\leq \psi_{l} \quad \forall u \in N, l \in \Omega \quad [-\pi_{ul}]\label{l_con_13}\\
    \sum_{u \in N}[d_u\geq d]y_{ul} \leq D_{dl}\psi_l\quad \forall l \in \Omega,d \in \mathcal{D} \quad [-\pi_{dl}]\label{l_con_23}\\
    \sum_{l \in \Omega}a_{fl}(\theta_l+\psi_l)\leq 1 \quad \forall f \in \mathcal{F} \quad [-\pi_f]\label{capDDOI3}
\end{align}
\end{subequations}
CG is used to solve \eqref{DT_DOI_master_primal3} in the same manner as in \eqref{DT_DOI_master_primal}.  Recall that pricing is performed only over $\theta_l$ which for SSCFLP is a knapsack problem as described in \eqref{SSCFLP_pricing}.  Recall that we never price over $y_{ul}$ or $\psi_l$.  This is because once $l$ is added to $\Omega_R$ via \eqref{SSCFLP_pricing} we add $\psi_l$ and $y_{ul} \; (\forall u \in N)$ to consideration in the RMP.

Observe that the DT-DOI provides the capacity that Flexible DOI (F-DOI) provide to the SSCFLP.  Specifically F-DOI provide the ability to recover the cost of servicing a customer $u$. This is interesting since DT-DOI are derived from Smooth-DOI, and not the F-DOI.  Note that the F-DOI are highly distinct from S-DOI as discussed in our literature review.

\section{Conclusion and Future Research}
\label{conc}

In this paper we introduced Detour-Dual Optimal Inequalities (DT-DOI) to stabilize the column generation (CG) approach to the Capacitated Vehicle Routing Problem (CVRP). DT-DOI provide the restricted master problem (RMP) the capacity to make detours along routes in the RMP to cover new items. The use of DT-DOI does not weaken the LP relaxation of the standard set cover formulation.
In future work we seek to incorporate DT-DOI into a branch-cut-price \citep{barnprice} framework. We also intend to use DT-DOI for problems with time windows, which is non-trivial since the insertion of detours may cause routes to violate time windows. However such a relaxation may be close to feasible permitting the use of relaxed DOI \citep{haghani2020relaxed} to ensure that an optimal and feasible solution to the CVRP with time windows is produced.  

We also intend to explore solving the DT-DOI RMP as an integer linear program to provide heuristic integer solutions at termination of CG for large problems where branch and cut and price can not be easily applied.  

Another topic for further research is that to adapt DT-DOI for the case where solving the RMP becomes a computational burden due to an explosion in the number of $y_{ul}$ terms over the course of CG.  In such a case we can decrease the number of $\psi_l$ terms, and $y_{ul}$ terms considered in the RMP.  The aim is to use only primal variables that we have reason to believe will be active in coming RMP solutions. For CVRP, one possibility is using all $\psi_l$ terms (for $l \in \Omega_R$)  but only the $y_{ul}$ terms for which $c_{ul}$ is small.  This is because terms with high $c_{ul}$ seem unlikely to be used when minimizing the RMP objective.  This approach can be adapted to Single Source Capacitated Facility Location (SSCFLP) where we also include $y_{ul}$ terms for which $a_{ul}=1$ so as to permit the original column to be expressed.  For both CVRP and SSCFLP we can choose to remove from consideration $\psi_{l}$ terms that have not been active in any solution to the RMP for several recent iterations since this would be a good indication that they will not be used in subsequent iterations.  Observe that the removal of $\psi_l$ from consideration removes all $y_{ul}$ from consideration.  We intend to explore these and other possibilities in future work.


\clearpage

\bibliographystyle{abbrvnat} 
\singlespacing
\bibliography{col_gen_bib}
\clearpage
\appendix

\section{Additional Figures}
\label{addFig}
For completeness, we present some additional figures from our experimental analysis. The results are similar to those presented earlier in Figs \ref{fig:iter1} and \ref{fig:iter2}.

\begin{figure}[!hbtp]
	\includegraphics[width=0.49\linewidth]{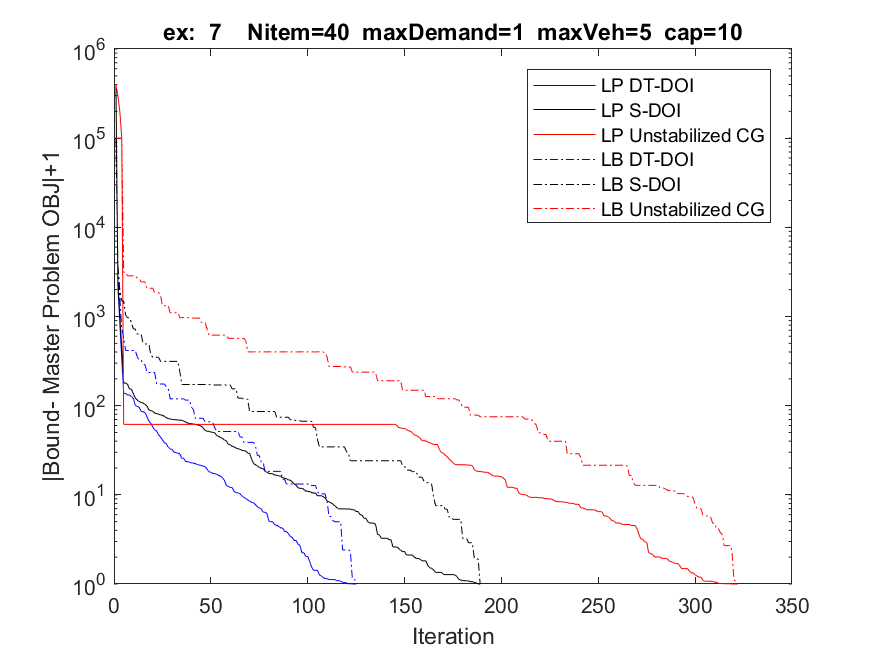}
	\includegraphics[width=0.49\linewidth]{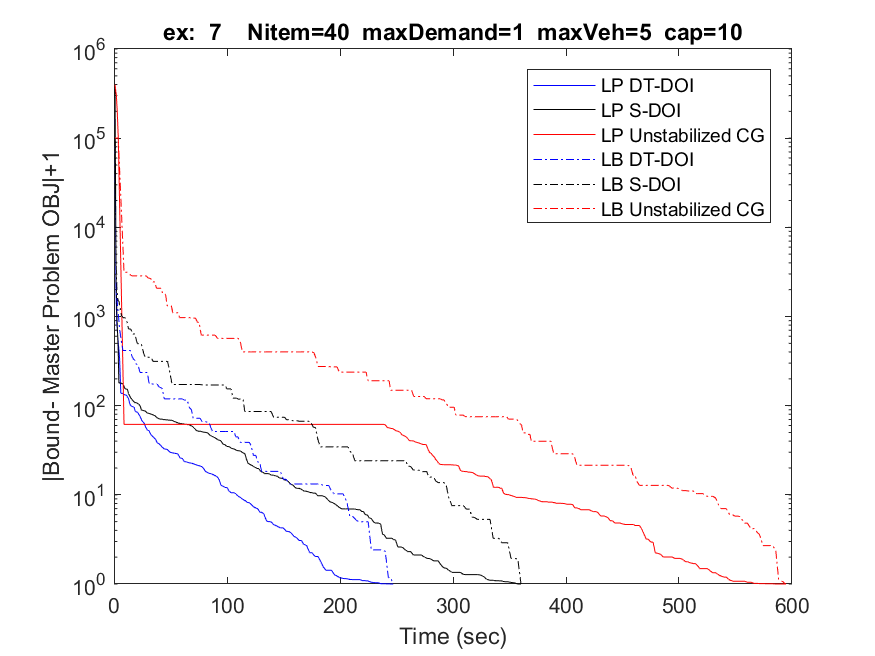}\\
	\includegraphics[width=0.49\linewidth]{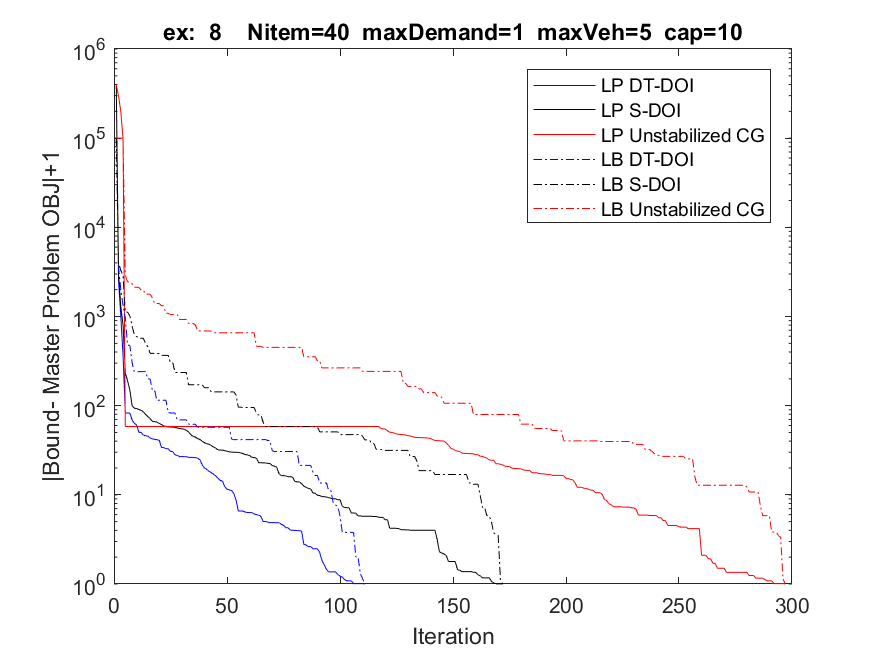}
	\includegraphics[width=0.49\linewidth]{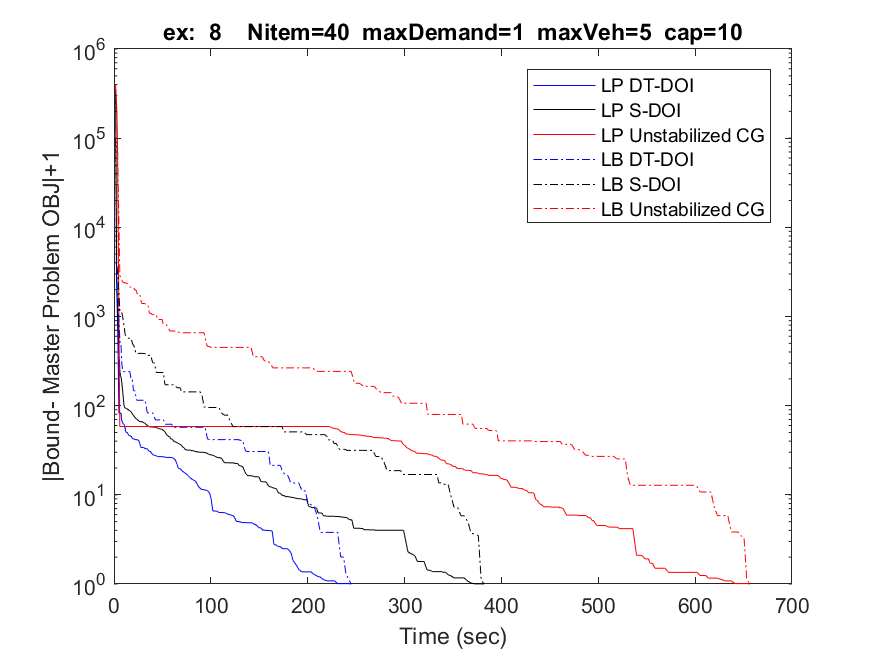}\\
	\includegraphics[width=0.49\linewidth]{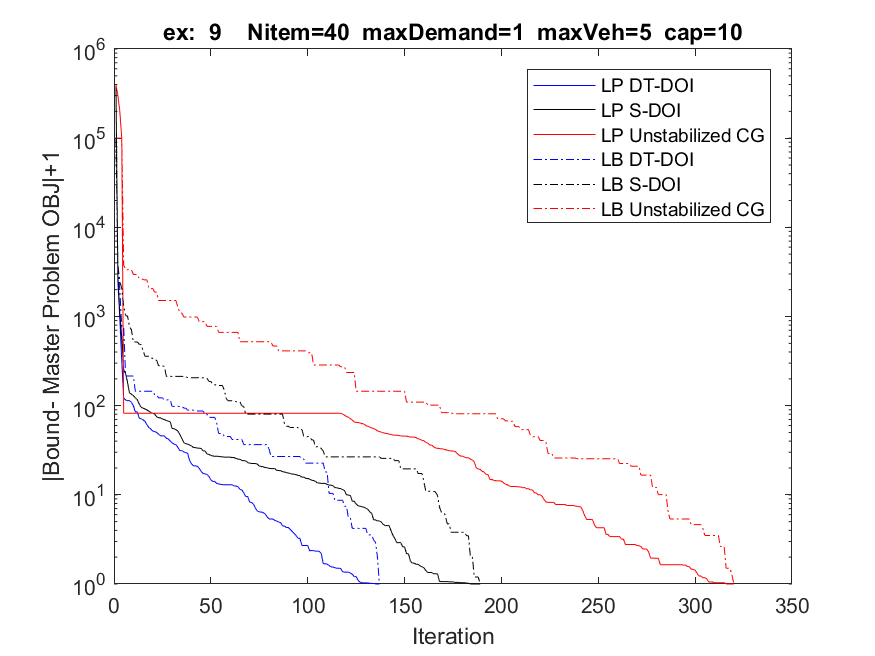}
	\includegraphics[width=0.49\linewidth]{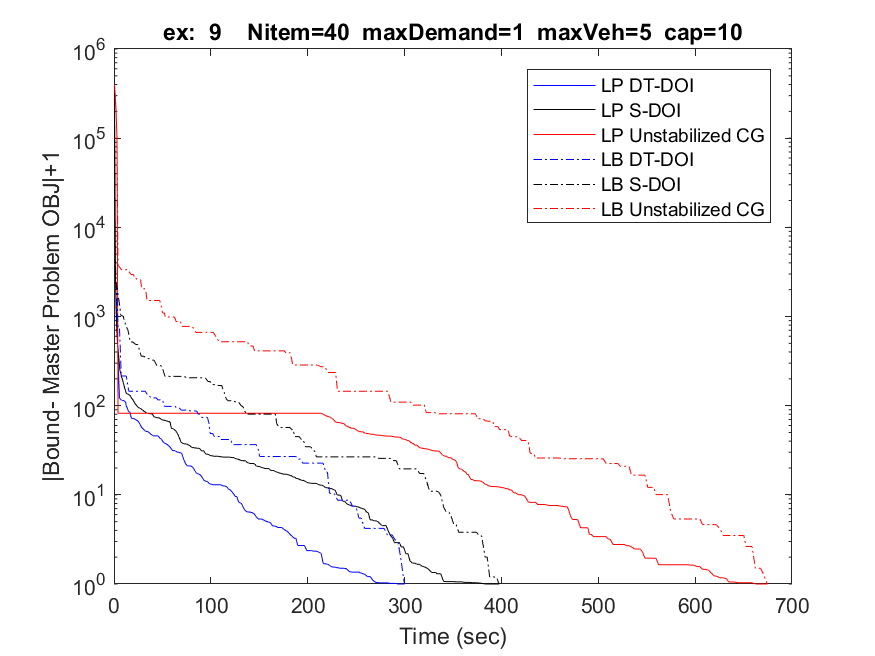}\\
	\caption{Results as a function of iteration/time.  The left side provides the results with respect to iteration and the right side for time on the same instance.  We add one to the difference of all bounds and MP values  which allows us to use the semilog scale.  
	}
	\label{fig:iter3}
\end{figure}
\begin{figure}[!hbtp]
	\includegraphics[width=0.49\linewidth]{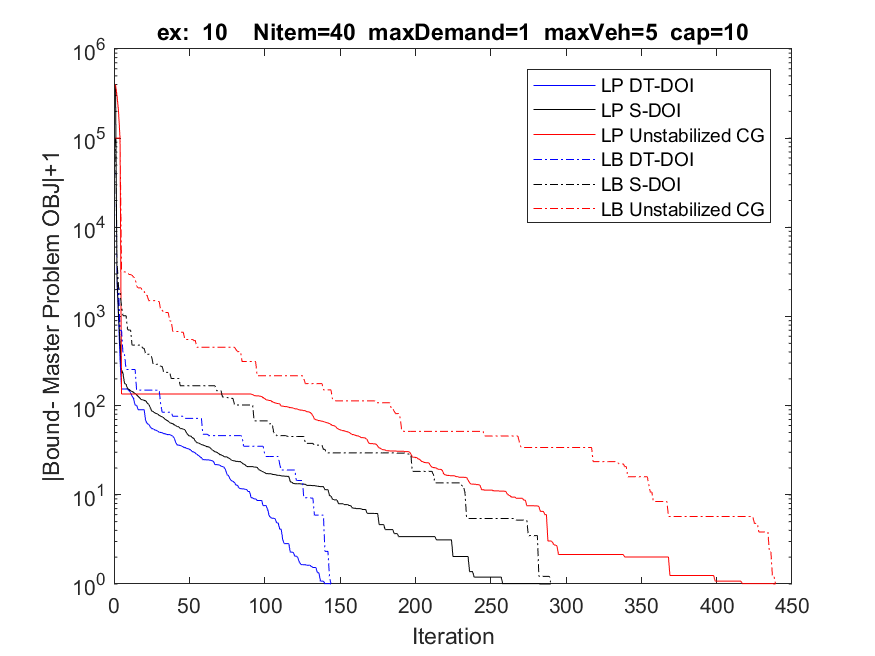}
	\includegraphics[width=0.49\linewidth]{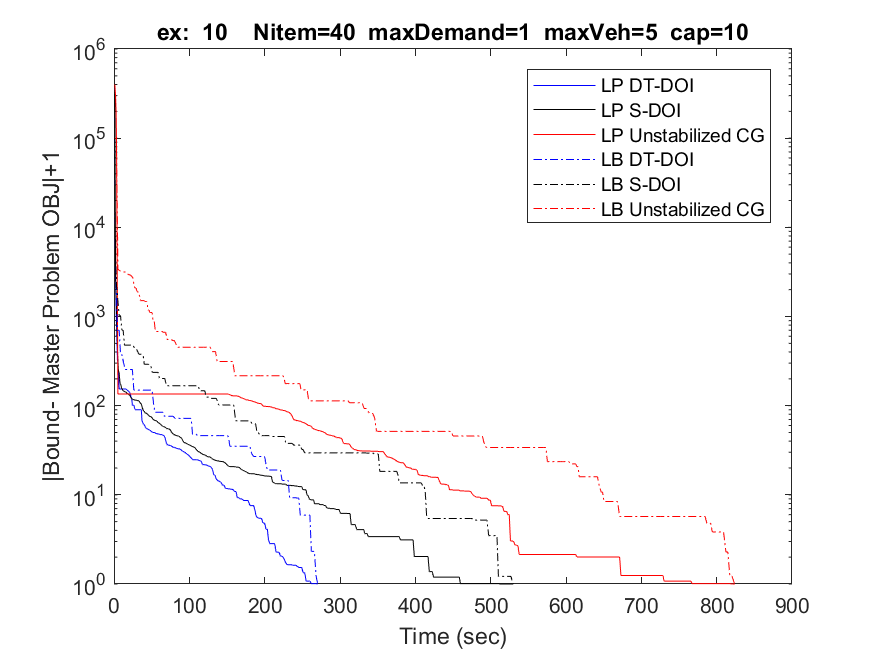}\\
	\includegraphics[width=0.49\linewidth]{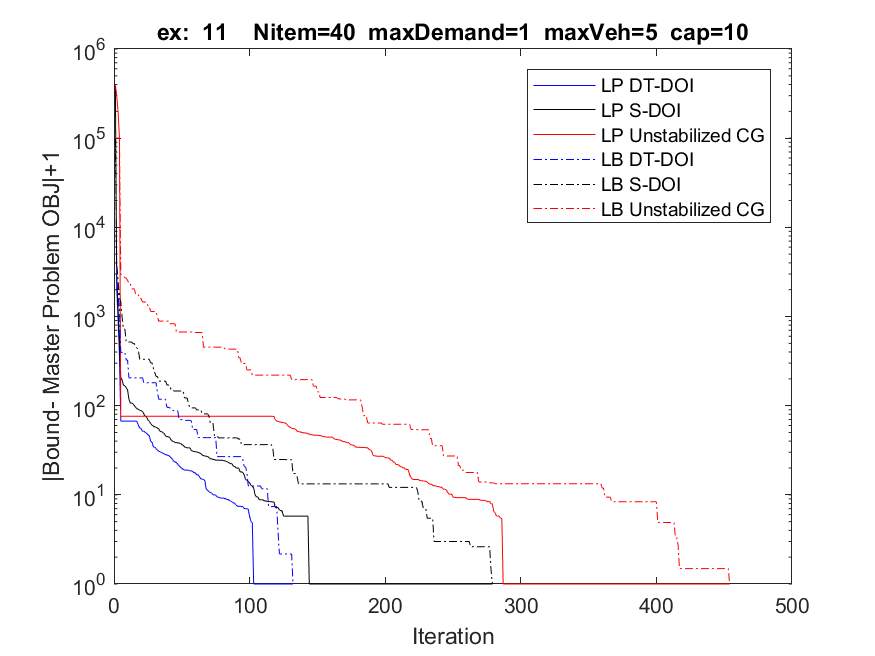}
	\includegraphics[width=0.49\linewidth]{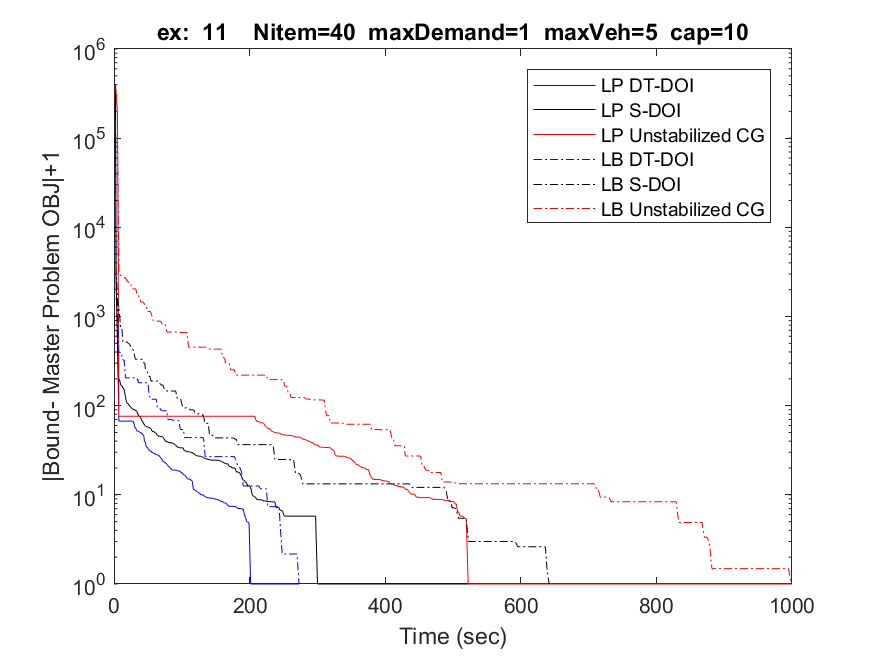}\\
	\includegraphics[width=0.49\linewidth]{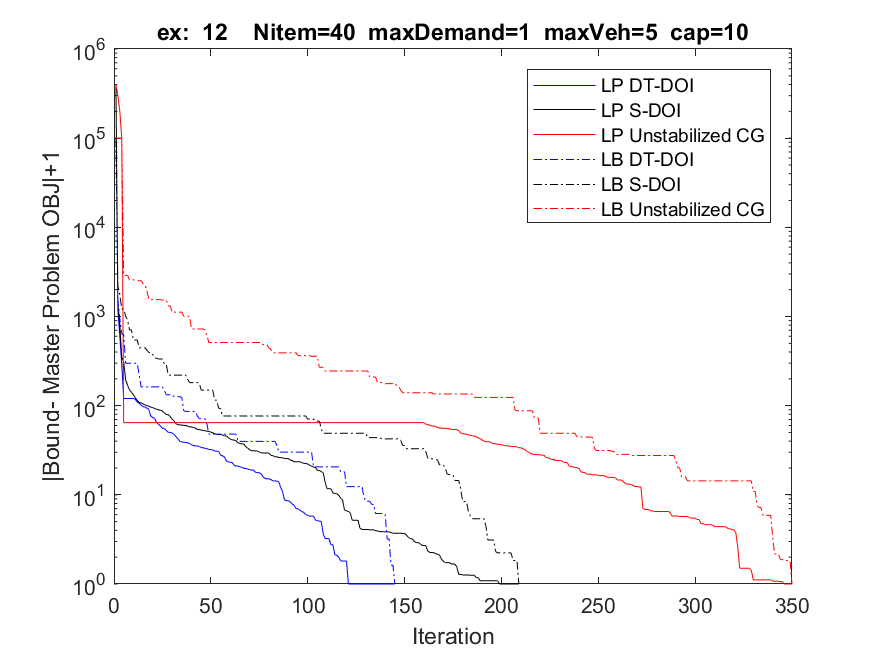}
	\includegraphics[width=0.49\linewidth]{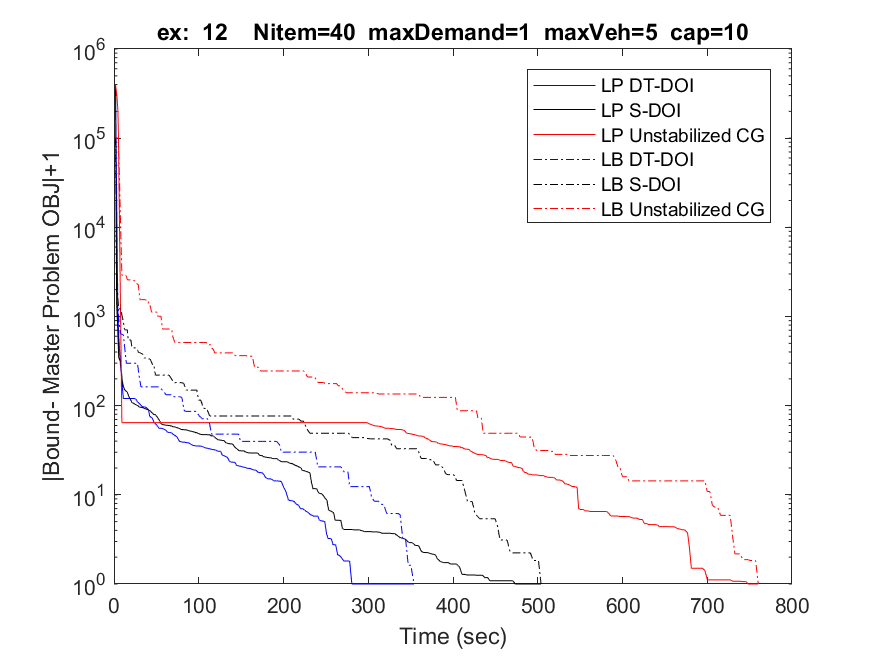}\\
	\caption{Results as a function of iteration/time.  The left side provides the results with respect to iteration and the right side for time on the same instance.  We add one to the difference of all bounds and MP values  which allows us to use the semilog scale.  
	}
	\label{fig:iter4}
\end{figure}
\begin{figure}[!hbtp]
	\includegraphics[width=0.49\linewidth]{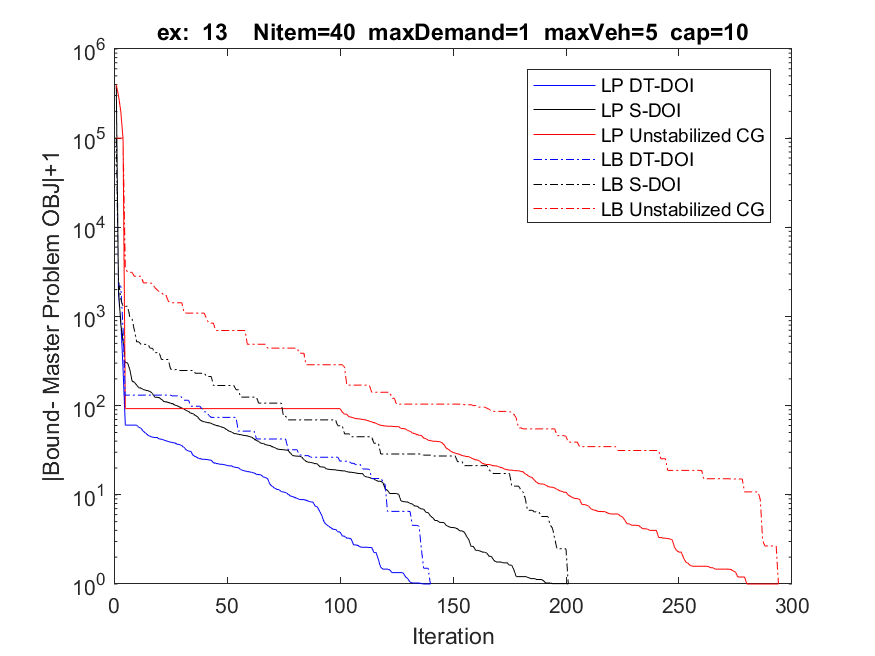}
	\includegraphics[width=0.49\linewidth]{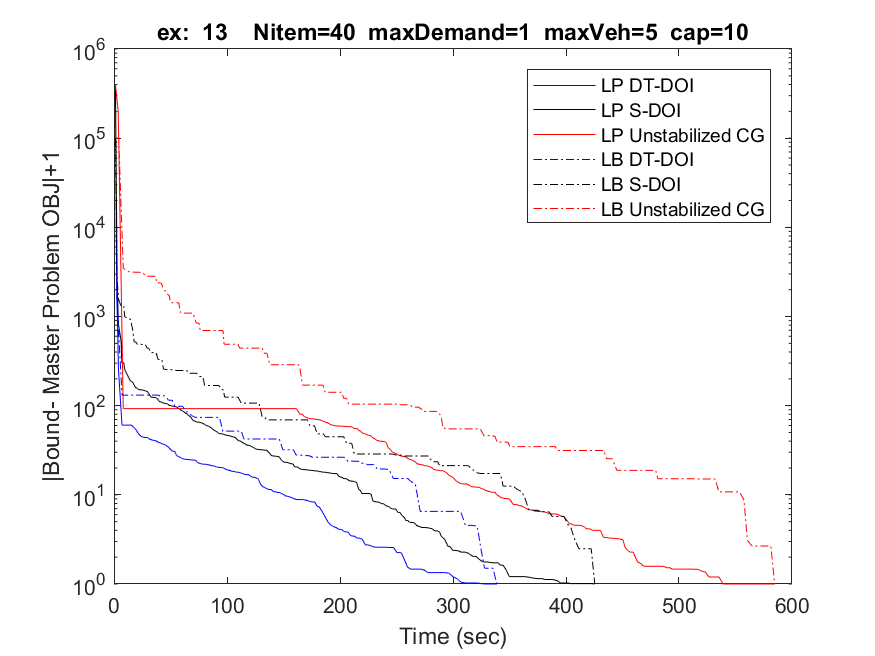}\\
	\includegraphics[width=0.49\linewidth]{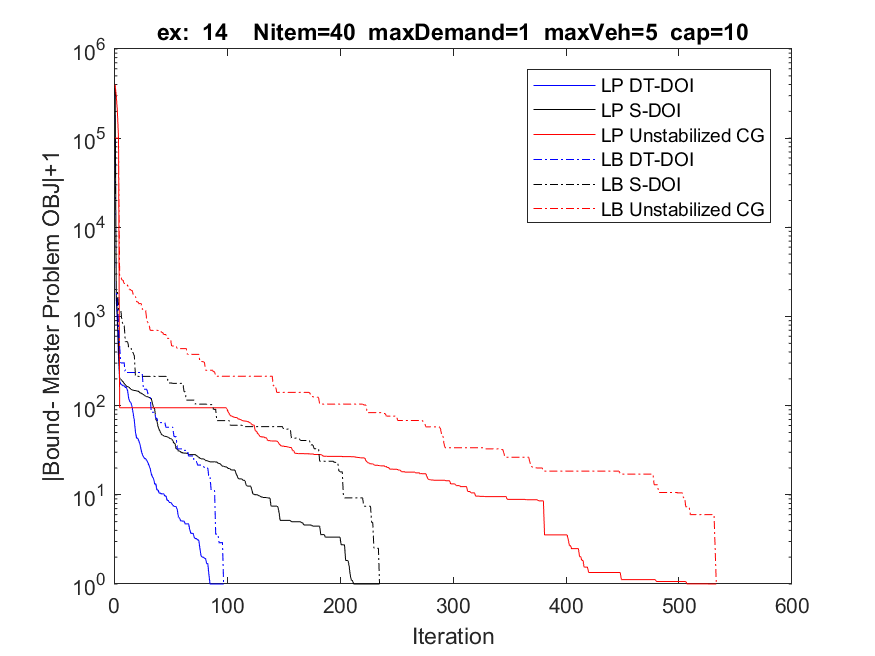}
	\includegraphics[width=0.49\linewidth]{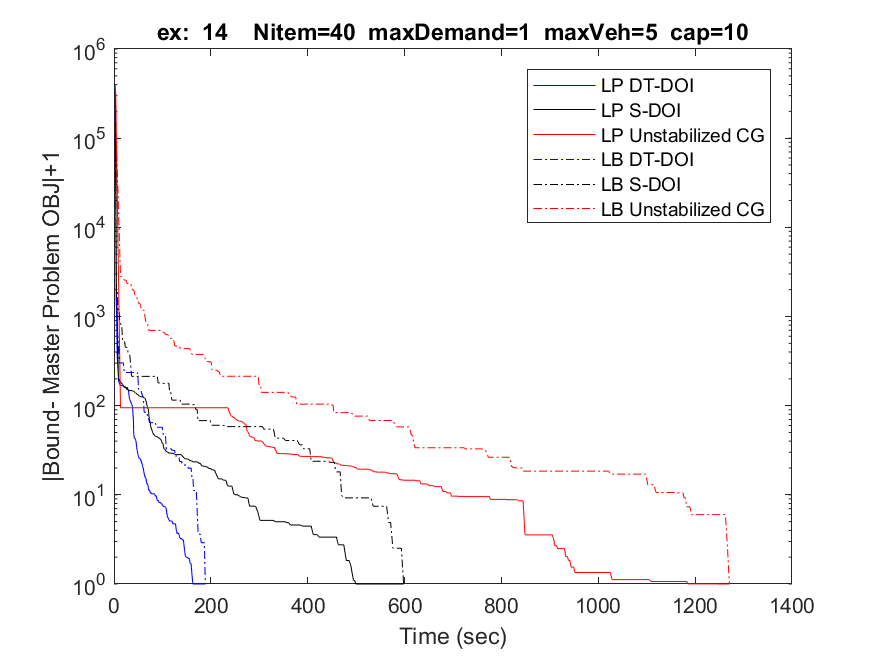}\\
	\caption{Results as a function of iteration/time.  The left side provides the results with respect to iteration and the right side for time on the same instance.  We add one to the difference of all bounds and MP values  which allows us to use the semilog scale.  
	}
	\label{fig:iter5}
\end{figure}

\section{Pricing as an Integer Linear Program}
\label{pricingILP}
We now consider the solution to pricing \eqref{pricing} as an integer linear program (ILP).

  We use decision variable $x_{uvd}=1$ if the generated route services $u$ then $v$ and contains (has remaining) exactly $d$ units of demand after leaving $u$ and otherwise set $x_{uvd}=0$. The following combinations of $u,v,d$ exist
\begin{itemize}
    \item $x_{uvd}$ exists if $d_v\leq d\leq D_0-d_u$ where $d_{-1}=d_{-2}=0$
\end{itemize}
The set of valid combinations of $u,v,d$ is denoted $P$.  
Below we  define $\bar{c}_{uv}$ to be the cost of traveling from $u$ to $v$ minus the additional cost corresponding to dual variables.  
\begin{subequations}
\begin{align}
    \bar{c}_{uv}=c_{uv}-\pi_v \quad  \forall (u,v,d) \in P, v \neq -2\\
    \bar{c}_{uv}=c_{uv}+\pi_0 \quad \forall (u,v,d) \in P, v = -2 \quad (\mbox{recall $v$=-2 is the end depot}) 
\end{align}
\end{subequations}
We write the lowest reduced cost resource constrained shortest path as pricing below in the form of an ILP, which we annotate after the ILP.  
\begin{subequations}
\label{pricingEqFull}
\begin{align}
    \min_{x \in \{0,1\} }\sum_{\substack{u,v,d \in P}}\bar{c}_{uv}x_{uvd} \label{ObjPricer}\\
    \sum_{\substack{u,v,d \in P}}x_{uvd}\leq 1 \quad \forall u \in N \label{elemnt}\\
    \sum_{\substack{-1,v,d \in P}}x_{-1vd}\leq 1 \label{pushCon}\\ 
    \sum_{\substack{u \in N}}x_{uvd}=\sum_{\substack{u \in N}}x_{v,u,d-d_v} \quad \forall  v\in N,d\geq d_v \label{flowcon}
\end{align}
\end{subequations}

In \eqref{ObjPricer} we minimize the reduced cost of the route.  In \eqref{elemnt} we ensure that an item is visited no more than once. In \eqref{pushCon} we ensure that up to no more than one route is selected by ensuring that no more than one unit of flow leaves the start depot.  In \eqref{flowcon} we ensure that the solution describes a route by enforcing the flow constraint.  

We should note that the solution of \eqref{pricingEqFull} is not typically solved as an ILP but instead tackled as with labeling algorithm \citep{costa2019} for sake of efficiency.
\section{Proof: DT-DOI does not Loosen the Relaxation }
\label{sec_dualProof}

In this section we establish that $\eqref{DT_DOI_master_primal}=\eqref{primal_master}$.  Clearly \eqref{DT_DOI_master_primal} is a relaxation of \eqref{primal_master} since any feasible solution to \eqref{primal_master} is a solution to \eqref{DT_DOI_master_primal} with the same objective.  We now establish that there is an optimal solution to \eqref{DT_DOI_master_primal} for which $\sum_{l \in \Omega}\psi_l=0$ and hence $y_{ul}=0 \quad \forall u \in N, l \in \Omega$, $\psi_l=0 \quad \forall l \in \Omega$. 

Consider the optimal solution to \eqref{DT_DOI_master_primal}.  If multiple optimal solutions exist select the solution that minimizes $\sum_{l \in \Omega}\psi_l$.  If $\psi=0 \quad \forall l \in \Omega$ then the claim is true since $y_{ul}$ is the zero vector by \eqref{l_con_1}. 
In this section we show that such an optimal solution to \eqref{DT_DOI_master_primal} for which $\sum_{l \in \Omega}\psi_l=0$ must exist using proof by contradiction.  In Section \ref{part1Proof} we show how to decrease $\sum_{l \in \Omega}\psi_l$ terms without increasing the objective \eqref{DDOIOBJ}  if  a specific condition is satisfied.  In Section \ref{part2Proof} we show that it this condition must  must be satisfied for any $l \in \Omega$ for which $\psi_l>0$.  Thus we establish a contradiction.
\subsection{Improving the Primal Solution}
\label{part1Proof}
In this section we consider any $l \in \Omega$ for which $\psi_l>0$.  We then alter the primal solution to decrease $\psi_l$,  without increasing \eqref{DDOIOBJ}, preserving feasibility and leaving $\psi_{\hat{l}}$ fixed for all $\hat{l}\in \Omega-l$.  Consider any $l$ for which $\psi_l>0$.  Let $v_{d},v_u$ be the slack on the constraint over $d,l$ in \eqref{l_con_2} and $u,l$ in \eqref{l_con_1} respectively. Let $\mathcal{D}^0,N^0$ be set of tight constraints over \eqref{l_con_2},\eqref{l_con_1} given $l$ respectively.  At least one constraint over \eqref{l_con_2},\eqref{l_con_1} must be tight otherwise $\psi_l$ could be decreased without altering feasibility or increasing the objective.  
 
Consider that there exists a $l_* \in \Omega$ satisfying the following.
\begin{subequations}
\label{DaRulesLstar}
\begin{align}
    a_{ul_*}=1 \quad \forall u \in N^0 \label{propallMustItem}\\
    a_{ul_*}=0 \quad \forall u \in N, y_{ul}=0 \label{prop_no_u}\\
     c_{l_*}\leq c_l+\sum_{u \in N_{l_*}}c_{ul}\label{propimprov}\\
     D_{dl_*}=D_{dl} \quad \forall d \in \mathcal{D}^0 \label{propSuff}
\end{align}
\end{subequations}
Now alter the solution as follows given step size $\alpha>0$ designed to ensure that $N^0,\mathcal{D}^0$ sets do not change and $y_{ul}\geq 0$ for all $u \in N$. We construct $\alpha$. 
\begin{subequations}
\label{myMod}
\begin{align}
\psi_{l}\leftarrow \psi_l-\alpha\\
\theta_{l_*}\leftarrow \theta_{l_*}+\alpha\\
y_{ul}\leftarrow y_{ul}-a_{ul_*}\alpha \quad \forall u \in N 
\end{align}
\end{subequations}
Let $\alpha$ be defined below as .999(any value $\in(0,1)$ works ) times the minimum of three terms positive valued terms $\alpha_1,\alpha_2,\alpha_3$ defined below with description subsequent.  
\begin{subequations}
\begin{align}
    \alpha=.999\min(\alpha_1,\alpha_{2},\alpha_3)\\
    \alpha_1=\min_{u \in N-N^0}v_u\\
    \alpha_2=\min_{d \in \mathcal{D}-\mathcal{D}^0}\frac{1}{D_{dl}}v_d\\
    \alpha_3=\min_{\substack{u \in N\\ y_{ul}>0}}y_{ul}
\end{align}
\end{subequations}
Here $\alpha_1$ is the minimum step size required to change $N^0$ when modifying according to \eqref{myMod}.  Similarly $\alpha_2$ is the minimum step size required to change $\mathcal{D}^0$ when modifying according to \eqref{myMod}.  We define $\alpha_3$ as a lower bound on the minimum step size required to change the set of $u$ for which $y_{ul}>0$.  Clearly $\alpha$ is positive since each of the component $\alpha_1,\alpha_2,\alpha_3$ are positive. 

Observe that the new solution describe in \eqref{myMod} remains primal feasible for \eqref{DT_DOI_master_primal} and has increased the objective by $\alpha(c_{l_*}- c_l-\sum_{u \in N_{l_*}}c_{ul})$ which is non-positive by \eqref{propimprov}.  We have also decreased $\psi_l$ by $\alpha$.  Thus we have  decreased the sum of the $\psi$ terms without increasing \eqref{DDOIOBJ} creating a contradiction with the claim that the unmodified solution is optimal and minimizes $\sum_{\hat{l} \in \Omega}\psi_{\hat{l}}$.  

\subsection{Establishing the Condition for Improving the Objective}
\label{part2Proof}

In this section we establish that for the $l$ in Section \ref{part1Proof} the $l_*$ as described in \eqref{DaRulesLstar} exists.  We first construct a set $\hat{N}$ whose elements will be identical to $N_{l_*}$ where $l_*$ is defined as the lowest cost route in $\Omega$ covering exactly the items in $\hat{N}$.

\textbf{Construction Procedure for $l_*$: }\textit{  We initialize $\hat{N}$ with all elements in $N^0$.  Now sort $u\in N-N^0$ for which $y_{ul}>0$ in order by $d_u$ from largest to smallest.  We iterate over these items.  When reaching an item $u$ we add it to $\hat{N}$ IFF $\hat{N}$ contains fewer than $D_{d_ul}$ items of demand $d_u$ or greater.  }
Note that we never add an item $u$ to $\hat{N}$ so that the number of items of demand $d_u$ or greater exceeds $D_{d_ul}$.  Thus all items in $\hat{N}$ can be serviced by a single route because $l$ is a feasible route.  Let $l_*$ be the lowest cost route servicing all items in $\hat{N}$.  

We now establish that $\hat{N}$ can be mapped to $l_*$ satisfying \eqref{DaRulesLstar}.  Since $\hat{N}$ is initialised to include all $u \in N^0$ then \eqref{propallMustItem} must hold.  Since no $y_{ul}$ for which $y_{ul}=0$ is ever added to $\hat{N}$ then \eqref{prop_no_u} holds.  

  Clearly \eqref{propimprov} holds since the RHS of \eqref{propimprov} describes a route that may visit additional items beyond $\hat{N}$; take detours; and or may not use the optimal ordering.  

By construction we know that \eqref{propSuff} is never violated by $D_{dl_*}>D_{dl}$ for any $d\in \mathcal{D}^0$.  We now establish  that $D_{dl_*}<D_{dl}$ for any $d \in \mathcal{D}^0$.  First observe that for any $d \in \mathcal{D}^0$ the following holds.
\begin{subequations}
\label{helper1Me}
\begin{align}
    \sum_{\substack{u \in N}}[d_u\geq d][y_{ul}>0]\psi_l
    \geq \sum_{\substack{u \in N}}[d_u\geq d][y_{ul}>0]\max_{\hat{u} \in N}y_{u\hat{l}}
    \geq \sum_{\substack{u \in N}}[d_u\geq d]y_{ul}=D_{dl}\psi_{l} \quad \forall d \in \mathcal{D}^0
\end{align}
\end{subequations}

Observe that  for any $d,d+1$ s.t.  $d \in \mathcal{D}^0$ and  $D_{dl}>D_{d+1,l}$ that the following holds.
\begin{subequations}
\label{helperMe2}
\begin{align}
    \sum_{\substack{u \in N}}([d_u\geq d]-[d_u\geq (d+1)])[y_{ul}>0]\psi_l\\
    \geq \sum_{\substack{u \in N}}([d_u\geq d]-[d_u\geq (d+1)])[y_{ul}>0]\max_{\hat{u} \in N}y_{u\hat{l}}\\
    \geq \sum_{\substack{u \in N}}([d_u\geq d]-[d_u\geq (d+1)])y_{ul}\\
    \geq (D_{dl}-D_{d+1,l})\psi_{l} \quad \forall d \in \mathcal{D}^0
\end{align}
\end{subequations}
Consider that during construction we satisfy all constraints of form \eqref{propSuff} up till $d$ for some $d\in \mathcal{D}$ where construction failed.  By fail we mean that a $u$ was reached for which $d_u<d$ prior to \eqref{propSuff} being satisfied for $d$.  We now show that such a failure can not occur.  

If $d$ is maximum demand in $\mathcal{D}$ then we know that \eqref{helper1Me} guarantees that there are at least $D_{dl}$ elements that can be added.  This is observed by dividing all terms in \eqref{helper1Me} by $\psi_l$ which is positive.  Thus the failure at $d$ can not occur.  

Failure can not occur at $d$ for which $D_{d+1,l}=D_{dl}$ is satisfied since adding zero items of demand $d$ satisfies the constraint in \eqref{propSuff} over $d$.

If failure occurs at $d$ for which $D_{d+1l}>D_{dl}$ then there are at least $D_{dl}-D_{d+1,l}$ elements in $N$ for which $([d_u\geq d]-[d_u\geq (d+1)])[y_{ul}>0]=1$; hence a failure can not occur at such a $d$.  We observe that there are at least $D_{dl}-D_{d+1,l}$ elements in $N$ for which $([d_u\geq d]-[d_u\geq (d+1)])[y_{ul}>0]=1$ by dividing \eqref{helperMe2} by $\psi_l$ which is positive.  

Since we have covered all cases for failure the construction algorithm has succeeded in producing an $l_*$ satisfying \eqref{DaRulesLstar}.   

\section{Decreasing the Size of DT-DOI Primal Problem}
\label{APP_removeThetaAppend}
 
We now consider an equivalent form of \eqref{DT_DOI_master_primal} that removes the $\theta$ terms without weakening the relaxation.  In experiments we employ this formulation where $\theta$ terms are ignored since it has few variables.  Any solution $\theta$ to \eqref{primal_master} can be mapped to one of the same cost using $\psi,y$ defined as follows.  $\psi_l \leftarrow \theta_l$ for all $l \in \Omega$; $y_{ul}\leftarrow a_{ul}\theta_l$ for all $u \in \mathcal{N},l \in \Omega$; then setting $\theta_l\leftarrow 0$ for all $l \in \Omega$.  .  We write this below using $\Omega_R$. 
\begin{subequations}
\label{DT_DOI_master_primal2}
\begin{align}
    \min_{\substack{\psi \geq 0 \\ y \geq 0}}\sum_{l \in \Omega_R}c_l\psi_l+\sum_{\substack{u \in N\\ l \in \Omega_R}}c_{ul}y_{ul} \label{DDOIOBJ2}\\
    \sum_{\substack{ l \in \Omega_R}}y_{ul} \geq 1 \quad \quad \forall u \in N \quad [\pi_u] \label{DDOIcover2}\\
    y_{ul}\leq \psi_{l} \quad \forall u \in N, l \in \Omega_R \quad [-\pi_{ul}]\label{l_con_12}\\
    \sum_{u \in N}[d_u\geq d]y_{ul} \leq D_{dl}\psi_l\quad \forall l \in \Omega_R,d \in \mathcal{D} \quad [-\pi_{dl}]\label{l_con_22}\\
    \sum_{l \in \Omega_R}\psi_l \leq K \quad [-\pi_0]\label{capDDOI2}
\end{align}
\end{subequations}
When solving optimization using \eqref{DT_DOI_master_primal2} as the RMP, pricing proceeds identically to using \eqref{DT_DOI_master_primal} (over $\Omega_R$) as the RMP.  Specifically we solve \eqref{pricing} generating $l$, which is then added to $\Omega_R$.  
\end{document}